\documentclass[10pt,a4paper,fleqn]{article}
\usepackage[latin1]{inputenc}
\usepackage{amsmath}
\usepackage{amsfonts}
\usepackage{amssymb}
\usepackage{color}
\usepackage[pdftex,plainpages=false,colorlinks,hyperindex,bookmarksopen,linkcolor=red,citecolor=blue,urlcolor=blue]{hyperref}
\usepackage[pdftex]{graphicx}

\usepackage[margin=4cm]{geometry}

\numberwithin{equation}{section}
\newtheorem{os}{Remark}[section]
\newtheorem{te}{Theorem}[section]

\renewcommand{\H}{H^{m,n}_{p,q}\left[ x \bigg| \begin{array}{l} (a_i, \alpha_i)_{i=1, .. , p}\\ (b_j, \beta_j)_{j=1, .. , q}  \end{array} \right]}

\author{Mirko D'Ovidio\\ Dipartimento di Statistica, Probabilit\`a e Statistica Applicata\\ Sapienza Universit\`a di Roma\\P.le Aldo Moro, 5\\00185 Rome (ITALY)\\ mirko.dovidio@uniroma1.it}

\date{July 28, 2009}

\title{Generalized Gamma Process: some results\\ about composition and subordination}

\begin{document}

\maketitle

\textbf{Abstract:} In this paper we deal with the generalized Gamma processes and their compositions. For the compositions of two or more than two generalized Gamma processes we give, when possible, the explicit law  whereas, in the other cases the representations in terms of Fox's H-functions are given. We also study the connections between iteration and product of random processes by exploiting the properties of the generalized Gamma processes, such a study allows us to obtain some striking result about the compositions of the Cauchy processes or fractional Brownian motions. Furthermore, we find out the partial differential equations governing the generalized Gamma processes and their compositions.\\


\tableofcontents

\section{Introduction}

A generalized Gamma random variable has a density law given by
\begin{equation}
\tilde{Q}(x;c,\mu,\gamma)= \gamma \frac{\left( \frac{x}{c} \right)^{\mu \gamma -1} e^{-\left( \frac{x}{c} \right)^{\gamma} }}{c \Gamma(\mu)}, \qquad  x>0 , \; c,\mu,\gamma >0
\label{gammaGen}
\end{equation}
where $c$ is a scale parameter and $\mu$, $\gamma$ are the shape parameters.\\
The generalized Gamma random variable is very important in the fields of reliability analysis and life testing models. It generalizes a number of distributions such as Weibull, Raleigh, folded normal, negative exponential. It is well-known that the Gamma random variable $G$ has a density law given by the formula \eqref{gammaGen} with $\gamma=1$, say $\tilde{Q}(x;c,\mu,1)$.\\ 
In literature, it refers to the Gamma process $\Gamma(t)$, $t>0$ that is governed by the density law 
\begin{equation}
q(x,t)=\frac{x^{t-1} e^{-x}}{\Gamma(t)}, \quad x>0,\; t>0.
\end{equation}
The process $\Gamma(t)$, $t>0$ is a subordinator, a L\'evy process with increasing paths and independent Gamma increments. Its density law is given by the formula \eqref{gammaGen} with the shape parameters $\mu=t$ and $\gamma=1$. In place of the process $\Gamma(t)$, $t>0$ we will consider, throughout the paper, the processes $\tilde{G}_\gamma(t)$, $t>0$ and $G_\gamma(t)$, $t>0$ which are slightly different from the process $\Gamma(t)$, $t>0$ discussed previously. The process $\tilde{G}_\gamma(t)$, $t>0$ has a density law given by the formula \eqref{gammaGen} with the shape parameter $c=t$. The distribution of the process $G_\gamma(t)$, $t>0$ is $q(x,t;\mu,\gamma)=\tilde{Q}(x; t^{1/\gamma}, \mu, \gamma)$ where $\tilde{Q}$ is that in \eqref{gammaGen}. However, all the processes in question are termed Gamma processes in the sense that $G_1(1)$ and $\Gamma(1)$ possess a Gamma distribution and $G_\gamma(1)$ and $\tilde{G}_\gamma(1)$ possess a generalized Gamma distribution.

We are interested in the study of the generalized Gamma processes $\tilde{G}_\gamma(t)$, $t>0$ and $G_\gamma(t)$, $t>0$ because they generalize a lot of well-known processes as the Brownian motion, the Bessel process starting from the origin, the exponential process and so on.
Furthermore, the compositions of generalized Gamma processes are connected to the Bessel functions, in particular, that of the second kind also known as the Macdonald function. It appears truly curious the role that the modified Bessel function $K_0(z)$ plays in the product of random processes.

We show that the distribution $q=q(x,t)$ of the Gamma process $G_1(t)$, $t>0$ solves the second order partial differential equation with non constant coefficients
\begin{equation}
\frac{\partial}{\partial t} q = x \frac{\partial^2}{\partial x^2}q - (\mu -2) \frac{\partial}{\partial x}q, \quad x>0, \; t >0. 
\end{equation}
Furthermore, for the generalized Gamma process $G_\gamma(t)$, $t>0$, we prove that the density law, say $q=q(x,t)$, solves the following partial differential equation
\begin{equation}
\frac{\partial}{\partial t} q = \frac{1}{\gamma^2} \left\lbrace \frac{\partial}{\partial x} x^{2-\gamma} \frac{\partial}{\partial x} q - (\gamma \mu -1) \frac{\partial}{\partial x} x^{1-\gamma} q \right\rbrace, \quad x >0, \; t>0
\label{eqUnoPDE}
\end{equation}
where $\gamma, \mu>0$. From the p.d.e. \eqref{eqUnoPDE} we obtain the partial differential equation governing the process $\tilde{G}_\gamma(t)$, $t>0$ which possesses distribution $Q(x;t^\gamma, \mu, \gamma)$. It suffices to consider that
\begin{equation}
\frac{\partial}{\partial t} Q(x;t^\gamma, \mu, \gamma) = \frac{\partial}{\partial z} Q(x;z, \mu, \gamma) \Bigg|_{z=t^\gamma} \frac{d (t^\gamma)}{dt}.
\end{equation}

In the present we study the composition of Gamma processes and we provide, when possible, an explicit form of their density law, otherwise representations in terms of the Fox's H-functions are given. For the composition $\tilde{G}^1_\gamma(\tilde{G}^2_\gamma(t))$, $t>0$ we prove that
\begin{equation}
\tilde{G}^1_\gamma(\tilde{G}^2_\gamma(t)) \stackrel{i.d.}{=} \tilde{G}^1_\gamma(t^{\frac{1}{2}}) \tilde{G}^2_\gamma(t^\frac{1}{2}), \quad t>0, \; \gamma \neq 0.
\end{equation}
Its density law is given by
\begin{equation}
q(x,t;\gamma,\mu)= \frac{2\gamma}{x \Gamma^2(\mu)} \left( \frac{x}{t} \right)^{\gamma \mu} K_0\left( 2 \sqrt{ \left( \frac{x}{t} \right)^\gamma} \right), \quad x >0, \; t>0,\; \mu >0, \; \gamma \neq 0
\label{fGGGG}
\end{equation}
where $K_0$ is the modified Bessel function with an imaginary argument.\\
We find out a very striking result about the composition $\tilde{G}^1_\gamma(\tilde{G}^2_{-\gamma}(t))$, $t>0$ which possesses the same distribution of the process $\tilde{G}^1_{-\gamma}(\tilde{G}^2_{\gamma}(t))$, $t>0$. The distribution in question reads
\begin{equation}
q(x,t;\mu,\gamma)= \gamma \frac{x^{\mu \gamma -1} t^{\mu \gamma}}{(x^\gamma + t^\gamma)^{2\mu}} \frac{\Gamma(2\mu)}{\Gamma^2(\mu)}, \quad x>0 \; t>0, \; \gamma >0.
\label{QaQa}
\end{equation}
Such processes generalize, in some sense, the Cauchy process $C(t)$, $t>0$ as it appears in \eqref{QaQa}. We also show that the following equivalences in distribution hold
\begin{equation}
\tilde{G}^1_\gamma(\tilde{G}^2_{-\gamma}(t)) \stackrel{i.d.}{=} \tilde{G}^1_\gamma(t^\frac{1}{2}) \tilde{G}^2_{-\gamma}(t^\frac{1}{2}), \quad t>0, \; \gamma \neq 0
\label{prUn}
\end{equation}
and
\begin{equation}
\tilde{G}^1_{-\gamma}(\tilde{G}^2_{\gamma}(t)) \stackrel{i.d.}{=} \tilde{G}^1_\gamma(t^\frac{1}{2}) \tilde{G}^2_{-\gamma}(t^\frac{1}{2}), \quad t>0, \; \gamma \neq 0.
\label{prDu}
\end{equation}
The last two equalities in distribution show the result mentioned before about the distribution \eqref{QaQa}. Indeed, the processes  \eqref{prUn} and \eqref{prDu} coincide and therefore they have the same distribution.\\
We also derive from the distribution \eqref{QaQa} some interesting result about the Student's t-distribution.\\
For $n$ independent generalized Gamma processes $\tilde{G}^j_{\gamma_j}(t)$, $t>0$, $j=1,2,\ldots , n$ we prove that
\begin{equation}
\tilde{G}^1_{\gamma_1}(\tilde{G}^2_{\gamma_2}(\ldots \tilde{G}^n_{\gamma_n}(t) \ldots)) \stackrel{i.d.}{=} \prod_{j=1}^n \tilde{G}^j_{\gamma_j}(t^\frac{1}{n}), \quad t>0, \; \gamma_j \neq 0, \, j=1,2,\ldots , n
\end{equation}
and the distribution can be expressed in terms of H-functions as follows
\begin{equation}
q(x,t)=\frac{(xt)^{-1}}{\Gamma^n(\mu)} H^{n,0}_{0,n} \left[ \frac{x}{t} \Bigg| \begin{array}{cccc} -;&-;& \ldots ;& -\\(\mu, \frac{1}{\gamma_1}) ;& (\mu, \frac{1}{\gamma_2}) ;& \ldots ;& (\mu, \frac{1}{\gamma_n}) \end{array} \right], \quad x \geq 0, \; t>0.
\end{equation}
In the special case where $\mu=\frac{1}{2}$ and $\gamma=2$, the function \eqref{QaQa} becomes the Cauchy density as we noted above. For the composition of the Cauchy process we prove that
\begin{equation}
a\,C^1(|C^2(\ldots |C^n(t)| \ldots)|) \stackrel{i.d.}{=} \prod_{j=1}^n  C^j \left( (at)^\frac{1}{n} \right) , \quad t>0, \; a>0
\end{equation}
where $C^j(t)$, $t>0$, $j=1,2,\ldots , n$ are independent Cauchy processes.\\

The process $G^1_{\gamma_1}(G^2_{\gamma_2}(t))$, $t>0$ has a representation of the density law in terms of Fox's H functions which reads
\begin{equation}
q(x,t) = \frac{1}{x\, t^{1/\gamma_1 \gamma_2}} H^{2,0}_{2,2}\left[ \frac{x}{ t^{1/\gamma_1 \gamma_2}} \Bigg| \begin{array}{ccc} (\mu, 0) & ; & (\mu , 0)\\ (\mu , \frac{1}{\gamma_1}) & ; & (\mu , \frac{1}{\gamma_1 \gamma_2}) \end{array} \right], \quad x>0,\; t>0
\end{equation}
where $\gamma_1, \gamma_2>0$. In the case where $\gamma_1=\gamma_2$ we have the distribution \eqref{fGGGG}  whereas in the special case where $\gamma_2=1$ the distribution has the explicit form 
\begin{equation}
q(x,t;\mu,\gamma_1)= 2\frac{\gamma_1 x^{\mu \gamma_1 -1}}{t^\mu \Gamma^2(\mu)} K_0\left( 2 \sqrt{ \frac{x^{\gamma_1}}{t} } \right), \quad x >0, \; t>0, \; \mu,\gamma_1 >0.
\label{iio}
\end{equation}
The iteration of $n$ independent Gamma processes $G^j_{\gamma_j}(t)$, $t>0$, $j=1,2\ldots, n$ leads to the process 
\begin{equation} 
I_G^{n}(t)= G^1_{\gamma_1} (G^2_{\gamma_2}( \ldots G^{n+1}_{\gamma_{n+1}}(t) \ldots)) , \quad t>0 , \; \gamma_j \neq 0 \; \forall j
\end{equation} 
which is related to some very appealing results.\\ 
We also show that the following equivalence in distribution holds
\begin{equation}
G^1_{\gamma_1} (G^2_{1}( \ldots G^{n}_{1}(t) \ldots)) \stackrel{i.d.}{=} \prod_{j=1}^{n} G_{\gamma_1}^j (t^{\frac{1}{n}}), \quad t>0
\label{dFh}
\end{equation}
where $|\gamma_1| >0$ and $\gamma_2=\gamma_3= \ldots = \gamma_n=1$. \\
The parameters $\gamma_1$ and $\mu_1$ turn out to be the most important parameters and one can ascribe such a matter to the fact that $G^1_{\gamma_1}(t)$, $t>0$ is the guiding process. We remind that $G^1_{\gamma_1}(t)$, $t>0$ posseses distribution $Q(x;t,\mu_1, \gamma_1)$.\\
The distribution of the process \eqref{dFh} can be expressed as
\begin{equation}
q(x,t;\mu)=\frac{x^{\gamma \mu-1}}{t^\mu \Gamma^n(\mu)} H^{n,0}_{0,n} \left[ \frac{x}{t^{1/\gamma}} \Bigg| \begin{array}{c} -\\ (0,\frac{1}{\gamma})_{i=1,2,\ldots , n} \end{array} \right], \quad x \geq 0, \; t>0
\end{equation}
where $\gamma_1=\gamma>0$ and $\mu_1=\mu>0$.

We also examine the process $B(G_1(t))$, $t>0$ where $B(t)$, $t>0$ is a standard Brownian motion and $G_1(t)$, $t>0$ is a Gamma process with distribution $Q(x;t,\mu,1)$. The process  $B(G_1(t))$, $t>0$ has the distribution 
\begin{equation}
q(x,t;\mu)= 2\frac{|x|^{\mu -\frac{1}{2}} }{ \pi^{\frac{1}{2}} \Gamma(\mu)} \left(\frac{2}{t} \right)^\frac{2\mu +1}{4} K_{\mu -\frac{1}{2}} \left( |x| \sqrt{\frac{2}{t}} \right), \quad x \in \mathbb{R},\; t>0, \; \mu>0.
\end{equation}
Such a function solves the p.d.e.
\begin{equation}
4 \frac{\partial}{\partial t}q = (2\mu -3) \frac{\partial^2}{\partial x^2} q - x \frac{\partial^3}{\partial x^3} q, \quad x \in \mathbb{R} \setminus \{ 0 \}, \; t>0.
\end{equation}
where we used the notation $q=q(x,t;\mu$) for the sake of simplicity.\\
It must be noted that compositions involving a fractional Brownian motion $B_H(t)$, $t>0$ and a generalized Gamma process $G_\gamma(t)$, $t>0$ turn out to have the same distribution of the process $B(G_1(t))$ for $\gamma = 2H$. In particular, we have
\begin{equation}
B_{H}(G_{2H}(t)) \stackrel{i.d.}{=} B(|G_{2H}(t)|^{2H}) \stackrel{i.d.}{=} B(G_1(t)), \quad t>0.
\end{equation}
Compositions involving the Gamma subordinator $\Gamma(t)$, $t>0$ are extensively studied in literature (see Madan, Carr, Chung \cite{MCC} and Madan, Seneta \cite{MS90} ). 

We study the $n-$dimensional process 
\begin{equation}
\left( B^1(G_1(t)), B^2(G_1(t)), \ldots , B^n(G_1(t)) \right)^T, \quad  t>0
\label{nDimBG}
\end{equation}
which is distributed as follows
\begin{equation}
q(\textbf{x},t) = \frac{2}{\pi^\frac{n}{2}} \frac{\| \textbf{x} \|^{\mu -\frac{n}{2}}}{\Gamma(\mu)}\left( \frac{2}{t}\right)^\frac{2\mu+n}{4} K_{\mu - \frac{n}{2}} \left( \| \textbf{x} \| \sqrt{\frac{2}{t}} \right), \quad \textbf{x} \in \mathbb{R}^n, \; t>0
\end{equation}
We give the p.d.e. governing the $n-$dimensional process in \eqref{nDimBG} which reads
\begin{equation}
4 \frac{\partial}{\partial t} q = (2\mu -n)\triangle q - \triangle ( \textbf{x} \cdot \nabla q ), \quad \textbf{x} \in \mathbb{R}^n, \; t>0 
\end{equation}
where $\triangle=\sum_{i=1}^n \partial_i^2 $ is the Laplacian operator, $\nabla=\sum_{i=1}^n \partial_i$ is the nabla operator and $q=q(\textbf{x},t)$.

Finally we study a general form of the Bessel process involving the generalized Gamma process. As already pointed out the generalized Gamma process $G_\gamma(t)$, $t>0$ possesses a density function which depends on the parameters $\mu$, $\gamma$ and $c=g(t)$ (for some positive $g$). This generalization allows to consider a large class of processes for which several results are proved to be true. We use the fact that
\begin{equation}
\left\lbrace \sum_{j=1}^{n} \left\lbrace  G^j_{\gamma_j,\mu}(t) \right\rbrace^{\gamma_j} \right\rbrace^\frac{1}{\gamma} \stackrel{i.d.}{=} G_{\gamma, n \mu}(t), \quad t>0.
\end{equation}

\section{Preliminaries}
We first introduce the Mellin transform of a sufficiently well-behaved function $f(x)$, $x \in (0,\infty)$. The Mellin transform of the function $f$ with parameter $\eta \in \mathbb{C}$ is defined as
\begin{equation}
\mathcal{M}\left\lbrace f(\cdot) \right\rbrace (\eta) = \int_0^\infty x^{\eta -1} f(x) dx
\end{equation}
whenever $\mathcal{M}\left\lbrace f(\cdot) \right\rbrace (\eta)$ exists. The inverse Mellin transform is defined as
\begin{equation}
f(x)=\frac{1}{2\pi i}\int_{ \theta -i \infty}^{\theta + i \infty} \mathcal{M}\left\lbrace f(x) \right\rbrace (\eta) x^{-\eta} d\eta 
\end{equation}
at all points $x$ where $f$ is continuous and for some real $\theta$.\\
Let us point out some useful operational rules that will be useful throughout the paper:
\begin{align}
& \mathcal{M}\left\lbrace f(a x) \right\rbrace (\eta) = a^{-\eta} \mathcal{M}\left\lbrace f \right\rbrace (\eta),\label{propM1}\\
& \mathcal{M}\left\lbrace x^a f(x) \right\rbrace (\eta) = \mathcal{M}\left\lbrace f \right\rbrace (\eta + a),\label{propM2}\\
& \mathcal{M} \left\lbrace f(ax^p) \right\rbrace (\eta) = \frac{1}{p} a^{-\frac{\eta}{p}} \mathcal{M}\left\lbrace f \right\rbrace \left(\frac{\eta}{p}\right), \\
& \mathcal{M}\left\lbrace \int_{0}^{\infty} f\left( \frac{x}{s} \right) g\left( s \right) \frac{ds}{s} \right\rbrace (\eta) = \mathcal{M}\left\lbrace f \right\rbrace (\eta) \cdot \mathcal{M}\left\lbrace g \right\rbrace (\eta) \label{conv}\\
& \mathcal{M} \left\lbrace \frac{\partial^n}{\partial x^n}f(x) \right\rbrace (\eta) = (-1)^n \frac{\Gamma(\eta)}{\Gamma(\eta - n)} \mathcal{M} \left\lbrace f \right\rbrace (\eta - n) \label{derM}.
\end{align}
The formula \eqref{conv} is the well-known Mellin convolution formula which turns out to be very useful in the study of the product of random variables.\\

Let us introduce the Fox's H-functions. The H-function is a very general function defined as follows
\begin{equation}
\H = \frac{1}{2\pi i} \int_{\mathcal{C}} \mathcal{H}^{m,n}_{p,q} (\eta) x^{\eta} d\eta
\end{equation}
where 
\begin{equation}
\mathcal{H}^{m,n}_{p,q} (\eta) = \frac{\prod_{j=1}^{m} \Gamma(b_j - \eta \beta_j) \prod_{i=1}^{n} \Gamma(1 - a_i + \eta \alpha_i)}{ \prod_{j_=m+1}^{q} \Gamma(1-b_j + \eta \beta_j) \prod_{i=n+1}^{p} \Gamma(a_j - \eta \alpha_i)}
\end{equation}
and $\mathcal{C}$ is a suitable path in the complex plane (for more details, see Mathai and Saxena \cite{MS78}).\\
Furthermore, the Mellin transform of a H-function reads
\begin{equation}
\int_{0}^{\infty} x^{\eta -1} \H dx = \mathcal{M}^{m,n}_{p,q} (\eta)
\end{equation}
where
\begin{equation}
\mathcal{M}^{m,n}_{p,q} (\eta)  = \frac{\prod_{j=1}^{m} \Gamma(b_j + \eta \beta_j) \prod_{i=1}^{n} \Gamma(1-a_i - \eta \alpha_i)}{\prod_{j=m+1}^{q} \Gamma(1-b_j - \eta \beta_j) \prod_{i=n+1}^p \Gamma(a_i + \eta \alpha_i)}.
\label{mellinHfox}
\end{equation}
In the sequel we need the following properties of the H-functions
\begin{equation}
\H = c\; H^{m,n}_{p,q}\left[ x^c \bigg| \begin{array}{l} (a_i, c \alpha_i)_{i=1, .. , p}\\ (b_j, c \beta_j)_{j=1, .. , q}  \end{array} \right] \quad c>0
\label{propH1}
\end{equation}
\begin{equation}
\H = \frac{1}{x^c} \; H^{m,n}_{p,q}\left[ x \bigg| \begin{array}{l} (a_i + c \alpha_i, \alpha_i)_{i=1, .. , p}\\ (b_j + c \beta_j, \beta_j)_{j=1, .. , q}  \end{array} \right], \quad c \in \mathbb{R}. \label{propH2}
\end{equation}

We also remind the modified Bessel function of the second kind $K_\nu$ which can be presented in several alternative forms: for $\Re\{ \nu \} > - \frac{1}{2}$ 
\begin{equation}
K_\nu (x)=\left( \frac{\pi}{2x} \right)^{\frac{1}{2}} \frac{e^{-x}}{\Gamma\left( \nu + \frac{1}{2} \right)} \int_{0}^{\infty} e^{-z} z^{\nu - \frac{1}{2}} \left( 1 + \frac{z}{2x} \right)^{\nu - \frac{1}{2}} dz, \quad | \arg x| < \pi 
\end{equation}
(see p. 140 Lebedev  \cite{LE}) or, for $\nu \neq 0, \; \nu=\pm 1, \pm2, \ldots$
\begin{equation}
K_\nu(x)=\frac{\pi}{2} \frac{I_{-\nu}(x) - I_{\nu}(x)}{\sin \nu \pi}, \quad | \arg x | < \pi
\end{equation}
(see p. 108 Lebedev  \cite{LE}) where
\begin{equation}
\label{I0Bessel1}
I_{\nu}(x)=\sum_{k=0}^{\infty} \left( \frac{x}{2} \right)^{2k + \nu} \frac{1}{k! \Gamma(k+\nu+1)}, \quad |x|<\infty, \; | \arg x | < \pi
\end{equation}
is the Bessel Modified function of the first kind (see p. 108 Lebedev \cite{LE}).\\ 
For $\nu \geq 0$ there exists the following representation
\begin{equation}
K_{\frac{\nu}{p}} \left(\sqrt{\frac{x}{t^{\zeta}}}\right) =\frac{| p |}{2} \left( x t^\zeta \right)^{-\frac{\nu}{2p}} \int_{0}^{\infty} s^{\nu -1} e^{-\frac{x}{2s^p} - \frac{s^p}{2 t^\zeta}}ds, \quad p \in \mathbb{R}, \; x>0, \; t>0, \; \zeta \in \mathbb{R}
\label{Knu}
\end{equation}
(see p. 370 Gradshteyn, Ryzhik \cite{GR}) and the special case where $\nu=0$, $p=2$
\begin{equation}
K_{0}(x)=\int_{0}^{\infty} s^{-1} e^{-\frac{x^2}{4s^2}-s^2} ds, \qquad | \arg x | < \frac{\pi}{4}
\label{functionK0}
\end{equation}
(see p. 119 Lebedev \cite{LE}). The Mellin transform of $K_\nu$ writes
\begin{equation}
\mathcal{M}\left\lbrace K_\nu \right\rbrace (\eta) = \frac{2^\eta}{4} \Gamma\left( \frac{\eta}{2} + \frac{\nu}{2} \right) \Gamma\left( \frac{\eta}{2} - \frac{\nu}{2} \right), \quad \Re\{\eta \} > \nu \geq 0. 
\label{wwwq}
\end{equation}
We observe that, for $\nu=0$, the Mellin transform \eqref{wwwq} becomes
\begin{equation}
\mathcal{M}\left\lbrace K_0  \right\rbrace (\eta) = \frac{2^\eta}{4} \Gamma^2\left( \frac{\eta}{2}\right), \quad \Re\{\eta \} > 0.
\end{equation}
Finally, in the following it will be useful the approximation
\begin{equation}
K_\nu(x) \approx \frac{2^{\nu -1} \Gamma(\nu)}{x^\nu}, \quad \textrm{for } x \to 0, \quad \nu >0
\label{K0in0}
\end{equation}
and
\begin{equation}
K_0 (x) \approx \log \frac{2}{x}, \quad x \to 0
\end{equation}
(see p. 136 Lebedev \cite{LE}).

\section{The generalized Gamma process}
In the literature, several authors worked with the Gamma process $\Gamma(t)$, $t>0$ which is a L\'evy process with increasing paths and c\`adl\`ag trajectories (a subordinator). Its density law is  given by
\begin{equation}
q(x,t)=\frac{x^{t-1} e^{-x}}{\Gamma(t)}, \quad x>0, \; t>0.
\label{gammaLevy}
\end{equation}
Hereafter, we deal with the generalized Gamma processes $\tilde{G}_\gamma(t)$, $t>0$ and $G_\gamma(t)$, $t>0$. The first process has a density law \eqref{gammaGen} and the second one is defined by the density
\begin{equation}
Q(x;c, \mu, \gamma) = \gamma \frac{x^{\mu \gamma -1} e^{-\frac{x^\gamma}{c}}}{c^\mu \Gamma(\mu)}, \quad x >0
\end{equation}
where $\gamma >0$, $\mu >0$ are the shape parameters and $c=g(t)$ is the scale parameter for some positive function $g:(0,\infty) \mapsto (0,\infty)$. The function $Q(x;c,\mu,\gamma)$ can be rewritten as
\begin{equation}
Q(x;c, \mu, \gamma) = \frac{1}{c^{1/\gamma}} Q_{\gamma, \mu} \left( \frac{x}{c^{1/\gamma}} \right)
\end{equation}
where, in explicit form
\begin{equation}
Q_{\gamma, \mu}(z) = \gamma \frac{z^{\mu \gamma -1} e^{-z^\gamma}}{\Gamma(\mu)}.
\end{equation}
The function \eqref{gammaGen} is derived as
\begin{equation}
Q(x;c^\gamma, \mu, \gamma) = \frac{1}{c}Q_{\mu, \gamma} \left( \frac{x}{c} \right)= \tilde{Q}(x;c,\mu,\gamma).
\end{equation}
\begin{os}
\normalfont
The function $q(x,t)=Q(x;t,\mu,\gamma)$ solves the first order partial differential equation
\begin{equation}
t \frac{\partial}{\partial t} q(x,t) = -\frac{1}{\gamma} \frac{\partial}{\partial x} \left( x q(x,t) \right), \quad x \in (0,\infty), \; t>0.
\end{equation}
Indeed, we have
\begin{equation*}
\frac{\partial}{\partial t} q(x,t) = -\frac{\mu}{t}q(x,t) + \frac{x^\gamma}{t^2} q(x,t)
\end{equation*}
and
\begin{align*}
-\frac{1}{\gamma} \frac{\partial}{\partial x} \left(x q(x,t) \right) = & - \frac{1}{\gamma}q(x,t) - \frac{x}{\gamma} \left\lbrace \frac{\mu \gamma -1}{x} q(x,t) - \frac{\gamma}{x} \frac{x^\gamma}{t} q(x,t) \right\rbrace \\
= & -\mu q(x,t) + \frac{x^\gamma}{t} q(x,t).
\end{align*}
\end{os}

We note that for a process $G_\gamma(t)$, $t>0$ (equivalently denoted by $G_{\gamma, \mu}(t)$, $t>0$ in order to indicate a process whose density is $Q(x;t, \mu, \gamma)$) with density law $q(x,t)=Q(x;t,\mu,\gamma)$ it is straightforward to show that 
\begin{equation*}
[G_1(t)]^{1/\gamma} \stackrel{i.d.}{=} G_\gamma(t), \quad t>0, \; \gamma >0 \quad \textrm{and} \quad [G_\gamma(t)]^\gamma \stackrel{i.d.}{=} G_1(t), \quad t>0, \; \gamma>0.  
\end{equation*}
Furthermore, it can be shown that 
\begin{equation}
[G_\gamma(t) ]^{-1} \stackrel{i.d.}{=} G_{-\gamma}(t), \quad t>0, \; \gamma >0.
\label{invaaaa}
\end{equation} 
where $G_{-\gamma}(t)$, $t>0$ possesses the inverse Gamma density law $-Q(x;t,\mu,-\gamma)$, $\mu>0$, $\gamma>0$. It will be useful in the sequel to write down the Mellin transform of the functions $ \pm Q(x;t,\mu, \pm \gamma)$. We have
\begin{align}
\mathcal{M}\left\lbrace Q(\cdot \, ;t,\mu,\gamma) \right\rbrace (\eta) = & \frac{\Gamma\left( \frac{\eta-1}{\gamma} + \mu \right)}{\Gamma(\mu)} t^\frac{\eta -1}{\gamma}, \quad \Re\{\eta \} > 1-\mu\gamma, \; \gamma >0 \label{mellinGgamma}\\
= & E\left\lbrace G_\gamma(t) \right\rbrace^{\eta -1}  \nonumber
\end{align}
and
\begin{align}
\mathcal{M}\left\lbrace - Q(\cdot \, ;t,\mu,-\gamma) \right\rbrace (\eta) = & \frac{\Gamma\left( \frac{1-\eta}{\gamma} + \mu \right)}{\Gamma(\mu)}t^\frac{1-\eta}{\gamma},\quad \Re\{ \eta \} <\gamma \mu +1, \; \gamma >0 \label{mellinGamma2}\\
= &  E\left\lbrace G_{-\gamma}(t) \right\rbrace^{\eta -1}. \nonumber
\end{align}
The Mellin transform of the function \eqref{gammaGen} is given by
\begin{align}
\mathcal{M}\left\lbrace \tilde{Q}(\cdot \, ;t,\mu, \gamma) \right\rbrace (\eta) = & \mathcal{M}\left\lbrace Q(x;t^{\gamma},\mu, \gamma) \right\rbrace (\eta) \label{melGGen}\\
= & \frac{\Gamma\left( \frac{\eta -1}{\gamma} + \mu \right)}{\Gamma(\mu)} t^{\eta -1}, \qquad \Re\{\eta \} > 1-\mu\gamma, \; \gamma>0 \nonumber \\
= & E\left\lbrace \tilde{G}_\gamma(t) \right\rbrace^{\eta-1} \nonumber
\end{align}
and, for completeness
\begin{align}
\mathcal{M}\left\lbrace - \tilde{Q}(\cdot \,;t,\mu, -\gamma) \right\rbrace (\eta) = & \mathcal{M}\left\lbrace Q(x;t^{-\gamma},\mu, -\gamma) \right\rbrace (\eta) \label{mellinGgeninv} \\
= & \frac{\Gamma\left( \frac{1-\eta}{\gamma} + \mu \right)}{\Gamma(\mu)} t^{\eta -1}, \quad \Re\{\eta \} < \gamma \mu +1,\; \gamma >0 \nonumber \\
= & E\left\lbrace \tilde{G}_{-\gamma}(t) \right\rbrace^{\eta -1} \nonumber.
\end{align}

For the Gamma process $G_1(t)$, $t>0$ with density law $q(x,t)=Q(x;t,\mu, 1)$ we can state the following theorem
\begin{te}
The function
\begin{equation}
q(x,t)= \frac{x^{\mu -1} e^{-\frac{x}{t}}}{t^\mu \Gamma(\mu)}, \quad x \geq 0,\; t>0
\end{equation}
say $q=q(x,t)$, solves the second order partial differential equation
\begin{equation}
\frac{\partial}{\partial t} q = x\frac{\partial^2}{\partial x^2} q - (\mu-2)\frac{\partial}{\partial x} q, \quad x \geq 0,\; t>0.
\label{PDEm1}
\end{equation}
\end{te}
\paragraph{First proof:}
The time derivative is given by
\begin{equation*}
\frac{\partial}{\partial t} q = \left( -\frac{\mu}{t} + \frac{x}{t^2} \right) q.
\end{equation*}
The first derivative w.r. to $x$ is given by
\begin{equation*}
\frac{\partial}{\partial x} q = \left( \frac{\mu -1}{x} - \frac{1}{t} \right) q
\end{equation*}
and thus, the second derivative w.r. to the space reads
\begin{align*}
\frac{\partial^2}{\partial x^2}q= & \left\lbrace  - \frac{\mu -1}{x^2} + \frac{(\mu-1)^2}{x^2} - 2\frac{\mu-1}{xt} + \frac{1}{t^2} \right\rbrace q\\
= & \left\lbrace - \frac{\mu -1}{x^2} + \frac{(\mu-1)^2}{x^2} - \frac{\mu-1}{xt} +\frac{1}{xt} + \frac{1}{x} \left( -\frac{\mu}{t} + \frac{x}{t^2} \right)  \right\rbrace q\\
= & \left\lbrace  (\mu -2 )\frac{\mu -1}{x^2} - \frac{\mu-1}{xt} +\frac{1}{xt} \right\rbrace  q + \frac{1}{x} \frac{\partial}{\partial t} q \\
= & \left\lbrace  \frac{(\mu -2 )}{x}\left( \frac{\mu -1}{x} - \frac{1}{t} \right) \right\rbrace  q + \frac{1}{x} \frac{\partial}{\partial t} q \\
= & \frac{\mu -2}{x} \frac{\partial}{\partial x} q + \frac{1}{x} \frac{\partial}{\partial t} q.
\end{align*}
This concludes the first proof.
\paragraph{Second Proof:}
We obtain the same result by performing the Mellin transform of the function
\begin{equation*}
\Psi(x) = x \frac{\partial^2}{\partial x^2} q - (\mu-2) \frac{\partial}{\partial x} q.
\end{equation*}
First of all, by the formula \eqref{derM}, we note that
\begin{equation*}
\mathcal{M}\left\lbrace x \frac{\partial^2}{\partial x^2} f(x) \right\rbrace (\eta) = \eta (\eta -1) \mathcal{M} \left\lbrace f \right\rbrace (\eta -1)
\end{equation*}
and
\begin{equation*}
\mathcal{M}\left\lbrace \frac{\partial}{\partial x} f(x) \right\rbrace =-(\eta -1) \mathcal{M} \left\lbrace f \right\rbrace (\eta -1)
\end{equation*}
for a function $f \in C^2(0,\infty)$ such that $\mathcal{M}\{ f \}$ exists.\\
By combining the previous results one obtains
\begin{align*}
\mathcal{M}\left\lbrace \Psi \right\rbrace (\eta) = & (\eta +\mu -2) (\eta -1) \mathcal{M} \left\lbrace Q(\cdot \,;t,\mu,1) \right\rbrace (\eta -1) \\
= & [ \textrm{by } \eqref{mellinGgamma} ] = (\eta +\mu -2) (\eta -1) \frac{\Gamma\left( \eta + \mu -2 \right)}{\Gamma(\mu)} t^{\eta -2}\\
= & (\eta -1) \frac{\Gamma\left( \eta + \mu -1 \right)}{\Gamma(\mu)} t^{\eta -2}\\
= & \frac{\partial}{\partial t} \frac{\Gamma\left( \eta + \mu -1 \right)}{\Gamma(\mu)} t^{\eta -1}.
\end{align*}
This concludes the second proof. $\blacksquare$\\

The generalized Gamma process is a $\frac{1}{\gamma}-$self-similar process. We can easily check that
\begin{equation}
G_\gamma(at) \stackrel{i.d.}{=} a^\frac{1}{\gamma} G_\gamma(t), \quad t>0,\; a>0, \; |\gamma | >0.
\label{ppppA}
\end{equation}
Indeed, by the Mellin transform \eqref{mellinGgamma}, we have that
\begin{equation*}
E\left\lbrace G_\gamma(at) \right\rbrace^{\eta -1} = \frac{\Gamma\left( \frac{\eta-1}{\gamma} + \mu \right)}{\Gamma(\mu)} (at)^\frac{\eta -1}{\gamma} = a^\frac{\eta -1}{\gamma} E \left\lbrace G_\gamma(t) \right\rbrace^{\eta-1} =  E \left\lbrace a^\frac{1}{\gamma} G_\gamma(t) \right\rbrace^{\eta-1}.
\end{equation*}
for $\gamma >0$.\\ 
Moreover, by the formula \eqref{mellinGamma2}, we have that
\begin{equation*}
E\left\lbrace G_{-\gamma}(at) \right\rbrace^{\eta -1} = \frac{\Gamma\left( \frac{1- \eta}{\gamma} + \mu \right)}{\Gamma(\mu)} (at)^\frac{1-\eta }{\gamma} = a^\frac{1-\eta}{\gamma} E \left\lbrace G_{-\gamma}(t) \right\rbrace^{\eta-1} =  E \left\lbrace a^{-\frac{1}{\gamma}} G_{-\gamma}(t) \right\rbrace^{\eta-1}.
\end{equation*}
for $\gamma >0$.\\
Thus, the Mellin transforms of the distribution of both members in \eqref{ppppA} coincide for $|\gamma|>0$.\\

We give now the partial differential equation governing the generalized Gamma process $G_\gamma(t)$, $t>0$ which possesses the distribution $q(x,t)=Q(x;t,\mu,\gamma)$.
\begin{te}
The density law of the process $G_\gamma(t)$, $t>0$ 
\begin{equation}
q(x,t) = \frac{x^{\gamma \mu -1} e^{-\frac{x^\gamma}{t}}}{t^\mu \Gamma(\mu)}, \quad x \geq 0, \, t>0, \; \mu>0, \; \gamma > 0  
\label{densm}
\end{equation}
say $q=q(x,t)$, satisfies the following p.d.e.
\begin{equation}
\frac{\partial}{\partial t} q = \frac{1}{\gamma^2} \left\lbrace \frac{\partial}{\partial x} x^{2-\gamma} \frac{\partial}{\partial x} q - (\gamma \mu -1) \frac{\partial}{\partial x} x^{1-\gamma} q \right\rbrace, \quad x \geq 0, \; t>0
\end{equation}
\label{teCharming}
\end{te}
\paragraph{First proof:}
The Mellin transform of the function \eqref{densm} reads
\begin{equation}
\Psi_t(\eta) = [ \textrm{by } \eqref{mellinGgamma} ]= \Gamma\left( \frac{\eta -1}{\gamma} + \mu \right) \frac{t^\frac{\eta -1}{\gamma}}{\Gamma(\mu)}, \quad \Re\{ \eta \} > 1-\gamma \mu .
\label{Zyyy}
\end{equation}
We perform the time derivative of the formula \eqref{Zyyy} and we obtain
\begin{align*}
\frac{\partial}{\partial t} \Psi_t(\eta) = &\frac{\eta -1}{\gamma} \Gamma\left( \frac{\eta -1}{\gamma} + \mu \right) t^\frac{\eta -\gamma -1}{\gamma}\\
=  & \frac{\eta -1}{\gamma} \left( \frac{\eta - \gamma -1 + \gamma \mu}{\gamma} \right) \Gamma\left( \frac{\eta - \gamma -1}{\gamma} + \mu \right) t^\frac{\eta - \gamma -1}{\gamma}\\
= & \frac{1}{\gamma^2} (\eta -1) ( \eta - \gamma -1 + \gamma \mu ) \Psi_t(\eta - \gamma)\\
= & \frac{1}{\gamma^2} (\eta -1) ( \eta - \gamma) \Psi_t(\eta - \gamma) +  \frac{1}{\gamma^2} (\eta -1) ( \gamma \mu -1 ) \Psi_t(\eta - \gamma)\\
= & \frac{1}{\gamma^2} \mathcal{M}\left\lbrace \frac{\partial}{\partial x} x^{2-\gamma} \frac{\partial}{\partial x} q \right\rbrace (\eta) - \frac{(\gamma \mu -1)}{\gamma^2} \mathcal{M} \left\lbrace \frac{\partial}{\partial x} x^{1-\gamma} q \right\rbrace (\eta).
\end{align*}
The inverse Mellin transform, according to the properties \eqref{derM} and \eqref{propM2}, yields the claimed result.
\paragraph{Second proof:} 
We first evaluate the derivatives w.r. to $x$
\begin{align*}
\frac{\partial}{\partial x} x^{2-\gamma} \frac{\partial}{\partial x} q = & \frac{\partial}{\partial x} x^{2-\gamma} \left\lbrace \frac{\gamma \mu -1}{x}Q - \frac{\gamma}{x} \frac{x^\gamma}{t} q \right\rbrace \\
= & \frac{\partial}{\partial x} \left\lbrace (\gamma \mu -1) x^{1-\gamma} q - \gamma x q \right\rbrace \\
= & \frac{(1-\gamma)(\mu\gamma -1)}{x^\gamma} q + x^{1-\gamma} (\mu \gamma -1) \frac{\partial}{\partial x}q - \frac{\gamma}{t} q - \gamma \frac{x}{t} \frac{\partial}{\partial x} q 
\end{align*}
and
\begin{equation*}
\frac{\partial}{\partial x} x^{1-\gamma} q = \frac{(1-\gamma)}{x^\gamma}q + x^{1-\gamma} \frac{\partial}{\partial x} q.
\end{equation*}
By combining the space derivatives, we obtain
\begin{align*}
\frac{\partial}{\partial x} x^{2-\gamma} \frac{\partial}{\partial x} q - (\gamma \mu -1) \frac{\partial}{\partial x} x^{1-\gamma} q = & -\frac{\gamma}{t} q - \gamma \frac{x}{t} \frac{\partial}{\partial x} q \\
= & -\frac{\gamma}{t}q - \gamma \frac{x}{t} \left\lbrace \frac{\gamma \mu -1}{x}q - \frac{\gamma}{x} \frac{x^\gamma}{t} q \right\rbrace \\
= & \gamma^2 \left\lbrace -\frac{\mu}{t}q + \frac{x^\gamma}{t^2}q \right\rbrace .
\end{align*}
By observing that
\begin{equation*}
\frac{\partial}{\partial t} q = -\frac{\mu}{t}q + \frac{x^\gamma}{t^2}q
\end{equation*}
the proof is completed. $\blacksquare$\\

\begin{os}
\normalfont
Consider the Theorem \ref{teCharming}. For $\gamma=1$ we obtain the p.d.e. \eqref{PDEm1}. For $\gamma=2$ and $\mu=\frac{1}{2}$ we have the heat-type equation. Indeed, the function 
\begin{equation*}
Q=Q(x;t,\frac{1}{2},2)=\frac{e^{-\frac{x^2}{t}}}{\sqrt{\pi t}},\quad  |x|\geq 0,\; t>0 
\end{equation*}
satisfies the p.d.e
\begin{equation*}
\frac{\partial}{\partial t} Q = \frac{1}{4} \frac{\partial^2}{\partial x^2} Q, \quad x \in \mathbb{R}, \; t>0.
\end{equation*}
We observe that the function $\hat{Q}=Q(x; g(t), \mu, \gamma)$, for some positive function $g:(0,\infty) \mapsto (0,\infty)$, satisfies the p.d.e.
\begin{equation}
\frac{\partial}{\partial t} \hat{Q} = \frac{1}{\gamma^2} \left\lbrace \frac{\partial}{\partial x} x^{2-\gamma} \frac{\partial}{\partial x} \hat{Q} - (\gamma \mu -1) \frac{\partial}{\partial x} x^{1-\gamma} \hat{Q} \right\rbrace \frac{d g}{d t}, \quad x \geq 0, \; t>0.
\end{equation}
Furthermore, for $\gamma=2$, $\mu=\frac{d}{2}$ and $g(t)=2t$, $t>0$ we obtain the p.d.e. governing the $d-$dimensional Bessel process starting from zero. 
\end{os}

\subsection{The $n$-dimensional Generalized Gamma process}

The multivariate Gamma process 
\begin{equation}
\mathbf{G}_{\gamma, \mu}({\bf t}) = \left( G^1_{\gamma, \mu}(t_1), \ldots , G^n_{\gamma, \mu}(t_n) \right), \qquad t_1, \ldots , t_n >0
\end{equation}
possesses dirtibution 
\begin{equation}
q(\mathbf{x}, \mathbf{t}) = \gamma^n \frac{\prod_{i=1}^n x_i^{\gamma \mu-1} e^{-\sum_{i=1}^n \frac{x_i^\gamma}{\phi t_i}}}{\phi^{\mu(n-1)} \prod_{i=1}^n t_i^\mu} \sum_{k \geq 0} \left( \frac{\rho}{\phi^n} \prod_{i=1}^n \frac{x_i^\gamma}{t_i} \right)^k \frac{1}{k!} \frac{(\mu)_k}{\Gamma^{n}(\mu + k)}  
\label{dist:multivG}
\end{equation}
where $\mathbf{x},\, \mathbf{t} \in \mathbb{R}^n_+$ and $\phi=1-\rho$, $0 \leq \rho \leq 1$. The simbol $(\mu)_k$ denotes the quantity
\begin{equation*}
(\mu)_0=1, \quad (\mu)_k=\frac{\Gamma(\mu+k)}{\Gamma(\mu)}, \; k=1,2,\ldots.
\end{equation*}
As a direct check shows all the marginals of the distribution \eqref{dist:multivG} are Gamma distributions. For $\rho=0$  the distribution \eqref{dist:multivG} takes the form
\begin{align*}
q(\mathbf{x}, \mathbf{t}) = & \gamma^n \frac{\prod_{i=1}^n x_i^{\gamma \mu-1} e^{-\sum_{i=1}^n \frac{x_i^\gamma}{t_i}}}{\Gamma^n(\mu) \prod_{i=1}^n t_i^\mu}\\
= & \prod_{i=1}^n \gamma \frac{x_i^{\gamma \mu -1} e^{-\frac{x_i^\gamma}{t_i}}}{t_i^\mu \Gamma(\mu)}\\
= & \prod_{i=1}^n q(x_i, t_i), \quad \mathbf{x},\, \mathbf{t} \in \mathbb{R}^n_+.
\end{align*}
For $n=2$ the distribution \eqref{dist:multivG} can be written as
\begin{equation}
q(x,y;t,s)=\gamma^2 \frac{\phi^\mu}{\rho^\mu} \frac{e^{-\frac{x^\gamma}{\phi t}-\frac{y^\gamma}{\phi s}}}{xy\, \Gamma(\mu)} \left( \sqrt{\frac{\rho}{\phi^2} \frac{x^\gamma y^\gamma}{st}} \right)^{\mu+1} I_{\mu-1} \left( 2 \sqrt{\frac{\rho}{\phi^2} \frac{x^\gamma y^\gamma}{st}} \right)
\end{equation}
After some calculations, we have
\begin{align*}
E \left\lbrace G_{\gamma, \mu}(t_1) G_{\gamma , \mu}(t_2) \right\rbrace = & \phi^\mu \frac{(\phi^2 t_1 t_2)^\frac{1}{\gamma}}{\Gamma(\mu)} \sum_{k \geq 0} \frac{\rho^k}{k!} \frac{\Gamma^2\left( \mu + k + \frac{1}{\gamma} \right)}{\Gamma(\mu + k)}\\
= & \phi^\mu (\phi^2 t_1 t_2)^\frac{1}{\gamma} \frac{\Gamma^2\left(\mu + \frac{1}{\gamma} \right)}{\Gamma^2(\mu)} \sum_{k \geq 0} \frac{\rho^k}{k!} \frac{\left( \mu +\frac{1}{\gamma} \right)_k \, \left( \mu +\frac{1}{\gamma} \right)_k}{(\mu)_k}\\
= & \phi^\mu (\phi^2 t_1 t_2)^\frac{1}{\gamma} \frac{\Gamma^2\left(\mu + \frac{1}{\gamma} \right)}{\Gamma^2(\mu)} F\left( \mu +\frac{1}{\gamma} , \mu +\frac{1}{\gamma} ; \mu ; \rho \right)
\end{align*}
where  
\begin{align}
F(a,b;c;z) = & \sum_{k \geq 0} \frac{z^k}{k!} \frac{(a)_k\, (b)_k}{(c)_k}, \quad \Re\{c\}> \Re \{ b \}>0, \quad |z|<1 \\
= & \frac{\Gamma(c)}{\Gamma(b) \Gamma(c-b)} \int_0^1 u^{b-1} (1-u)^{c-b-1}(1-uz)^{-a}du
\end{align}
is the hypergeometric function (for more details, see p.238 \cite{LE}).\\
Also, for a genaralized Gamma process $G_{-\gamma, \mu}(t)$, $t>0$, we can write down
\begin{equation*}
E \left\lbrace G_{-\gamma, \mu}(t_1) G_{-\gamma , \mu}(t_2) \right\rbrace = ( t_1 t_2)^{-\frac{1}{\gamma}} \frac{\Gamma^2\left(\mu - \frac{1}{\gamma} \right)}{\Gamma^2(\mu)} F\left(\frac{1}{\gamma}, \frac{1}{\gamma}; \mu; \rho \right)
\end{equation*}
thanks to the property
\begin{equation*}
F(a,b;c;z) = (1-z)^{c-a-b} F(c-a, c-b; c; z)
\end{equation*}
(see A S p. 559 \cite{XX}).\\
Furthermore, for $\gamma=1$, it is immediate to show that
\begin{align*}
E \left\lbrace G_{1, \mu}(t_1) G_{1 , \mu}(t_2)\right\rbrace  = & \phi^\mu \frac{\phi^2 t_1 t_2}{\Gamma(\mu)} \sum_{k \geq 0} \frac{\rho^k}{k!} (\mu +k) \Gamma(\mu + k + 1)\\
= & \phi^\mu \frac{\phi^2 t_1 t_2}{\Gamma(\mu)} \frac{\mu}{\phi^{\mu+2}} (\phi \mu + \mu \rho + \rho)\\
= & \mu (\mu + \rho)t_1 t_2.
\end{align*}
The calculations above can be carried out by rewriting the gamma function as
\begin{equation*}
\Gamma(\mu+k+1) = \int_0^\infty u^{\mu+k} e^{-u}du.
\end{equation*}
By observing that $EG_{1,\mu}(t)=\mu t$, we obtain the covariance function of a Gamma process $G_{1,\mu}(t)$, $t>0$ given by
\begin{equation*}
R_1(t_1, t_2) = \mu \rho t_1 t_2, \qquad t_1,\, t_2 >0, \quad \mu>0 , \quad 0 \leq \rho < 1
\end{equation*}
which is only positive.

\begin{os}
\normalfont
We point out that a fractional Brownian motion has covariance function given by
\begin{equation*}
R_H (t_1, t_2) = \frac{\sigma^2}{2}\left( |t_1|^{2H} + |t_2|^{2H} - |t_1 - t_2|^{2H} \right)
\end{equation*}
which becomes $R_1(t_1, t_2) =\sigma^2 t_1 t_2$ in the special case where $H=1$.
\end{os}

\section{Compositions involving Generalized Gamma processes}

Compositions involving two generalized Gamma processes are now examined. Let us start from the composition of two independent Gamma processes $\tilde{G}^j_\gamma(t)$, $t>0$, $j=1,2$ each with density \eqref{gammaGen} and $c=t$. Our purpose is to derive an explicit form of the density of the compound process $\tilde{G}^1_\gamma (\tilde{G}^2_\gamma(t))$, $t>0$, that can be also viewed as $G^1_\gamma(|G^2_\gamma(t^\gamma)|^\gamma)$, $t>0$ where the compositions of the processes $G_\gamma(t)$, $t>0$ will be discussed later on.

\paragraph{The process $\tilde{G}^1_\gamma(\tilde{G}^2_\gamma (t))$.}

The density of the process $\tilde{G}^1_\gamma (\tilde{G}^2_\gamma(t))$, $t>0$ writes
\begin{align}
q(x,t;\gamma,\mu) = & \gamma^2 \frac{x^{\mu \gamma -1}}{t^{\mu\gamma} \Gamma^2(\mu)} \int_0^\infty s^{-1} e^{-\left( \frac{x}{s}\right)^\gamma} e^{-\left( \frac{s}{t} \right)^\gamma} ds \label{Gtilde} \\
= & [ \textrm{by } \eqref{Knu} ] \nonumber \\
= & \frac{2|\gamma |}{x \Gamma^2(\mu)} \left( \frac{x}{t} \right)^{\gamma \mu} K_0\left( 2 \sqrt{ \left( \frac{x}{t} \right)^\gamma} \right), \quad x >0, \; t>0,\;  \mu >0, \; |\gamma |>0 \nonumber
\end{align}
where $K_0$ is the modified Bessel function \eqref{functionK0}. The Mellin transform of the function \eqref{Gtilde} reads
\begin{equation}
\mathcal{M}\left\lbrace q(\cdot \,,t;\gamma,\mu) \right\rbrace (\eta) = E\left\lbrace \tilde{G}^1_\gamma (\tilde{G}^2_\gamma(t)) \right\rbrace^{\eta -1} = \frac{\Gamma^2\left( \frac{\eta -1}{\gamma} + \mu \right)}{\Gamma^2(\mu)} t^{\eta -1}, \quad \Re \{\eta \} >1-\gamma \mu.
\label{hhg}
\end{equation}
In light of the Mellin convolution formula \eqref{conv} we can argue that 
\begin{equation}
E\left\lbrace \tilde{G}^1_\gamma (\tilde{G}^2_\gamma(t)) \right\rbrace^{\eta -1} = E \left\lbrace  \tilde{G}^1_\gamma (t^\frac{1}{2}) \right\rbrace^{\eta -1} E\left\lbrace  \tilde{G}^2_\gamma(t^\frac{1}{2}) \right\rbrace^{\eta -1} 
\end{equation}
where $E\left\lbrace  \tilde{G}^i_\gamma(t^\frac{1}{2}) \right\rbrace^{\eta -1}$ for $i=1,2$ is given by the formula \eqref{melGGen} and thus, we can write down the following equality in distribution
\begin{equation}
\tilde{G}^1_\gamma (\tilde{G}^2_\gamma(t)) \stackrel{i.d.}{=} \tilde{G}^1_\gamma(t^\frac{1}{2}) \tilde{G}^2_\gamma(t^\frac{1}{2}), \quad t>0, \; |\gamma | >0.
\label{qBn}
\end{equation}
Thus, from the relation \eqref{qBn} we obtain  the following equalities
\begin{equation}
\tilde{G}^1_\gamma (\tilde{G}^2_\gamma(t)) \stackrel{i.d.}{=} \tilde{G}^1_\gamma(t^\frac{1}{2}) \tilde{G}^2_\gamma(t^\frac{1}{2}), \quad t>0, \; \gamma >0
\end{equation}
and
\begin{equation}
\tilde{G}^1_{-\gamma} (\tilde{G}^2_{-\gamma}(t)) \stackrel{i.d.}{=} \tilde{G}^1_{-\gamma}(t^\frac{1}{2}) \tilde{G}^2_{-\gamma}(t^\frac{1}{2}), \quad t>0, \; \gamma >0.
\label{sFg}
\end{equation}
The distribution of the processes in \eqref{sFg} is given by
\begin{equation}
q(x,t;-\gamma, \mu) = \frac{2\gamma }{x \Gamma^2(\mu)} \left( \frac{t}{x} \right)^{\gamma \mu} K_0 \left( 2\sqrt{\left( \frac{t}{x} \right)^\gamma} \right), \quad x >0, \; t>0, \; \gamma >0
\end{equation}
which is derived from the formula \eqref{Gtilde}.

\paragraph{The process $\tilde{G}^1_\gamma(\tilde{G}^2_{-\gamma}(t))$.}

We consider now the compound process $\tilde{G}^1_\gamma(\tilde{G}^2_{-\gamma}(t))$, $t>0$ obtained by composing two independent Gamma processes. Such processes are, respectively, a generalized Gamma process $\tilde{G}_\gamma(t)$, $t>0$ and an inverse generalized Gamma process $\tilde{G}_{-\gamma}(t)$, $t>0$. The density law of the process $\tilde{G}^1_\gamma(\tilde{G}^2_{-\gamma}(t))$, $t>0$ writes
\begin{align}
q(x,t;\mu,\gamma) = & \gamma^2 \int_0^\infty \left( \frac{x}{s} \right)^{\mu \gamma -1} \frac{e^{-\left( \frac{x}{s} \right)^\gamma}}{s \Gamma(\mu)} \left( \frac{s}{t} \right)^{- \mu \gamma -1} \frac{e^{-\left( \frac{s}{t} \right)^{-\gamma}}}{t \Gamma(\mu)} ds \label{aaw} \\
= & \gamma^2 \frac{(xt)^{\gamma \mu}}{x \Gamma^2(\mu)} \int_{0}^{\infty} s^{-2\gamma \mu -1} e^{-\left( \frac{x}{s} \right)^\gamma - \left( \frac{s}{t} \right)^{-\gamma}} ds  \nonumber \\
= & \gamma^2 \frac{(xt)^{\gamma \mu}}{x \Gamma^2(\mu)} \int_{0}^{\infty} s^{2\gamma \mu -1} e^{-s^\gamma \left( x^\gamma + t^\gamma \right)} ds \nonumber \\
= & \gamma \frac{x^{\mu \gamma -1} t^{\mu \gamma}}{(x^\gamma + t^\gamma)^{2\mu}} \frac{\Gamma(2\mu)}{\Gamma^2(\mu)}, \quad x >0, \; t>0, \; \mu,\gamma>0. \nonumber
\end{align}
The Mellin transform of the function \eqref{aaw} reads
\begin{equation}
\mathcal{M}\left\lbrace q(\cdot \,,t;\mu,\gamma) \right\rbrace (\eta) = \Gamma\left(\frac{\eta-1}{\gamma} +\mu \right) \Gamma\left(\frac{1-\eta}{\gamma} + \mu \right) \frac{t^{\eta-1}}{\Gamma^2(\mu)}, \quad 0< \Re\{\eta\} < \gamma \left( \mu +\frac{1}{\gamma} \right)
\label{mellinaaw}
\end{equation}
and then, we can ask ourselves if the Mellin convolution results still hold.  It is straightforward to rewrite  the transform \eqref{mellinaaw} as
\begin{equation}
E \left\lbrace \tilde{G}^1_\gamma(\tilde{G}^2_{-\gamma}(t)) \right\rbrace^{\eta -1} = E\left\lbrace \tilde{G}^1_\gamma(t^\frac{1}{2}) \right\rbrace^{\eta-1}  E \left\lbrace  \tilde{G}^2_{-\gamma}(t^\frac{1}{2}) \right\rbrace^{\eta -1} .
\label{mellinQQQQ}
\end{equation}
Therefore, the following equality in distribution holds
\begin{equation}
\tilde{G}^1_\gamma(\tilde{G}^2_{-\gamma}(t)) \stackrel{i.d.}{=} \tilde{G}^1_\gamma(t^\frac{1}{2}) \tilde{G}^2_{-\gamma}(t^\frac{1}{2}), \quad t>0.
\label{productGtilde}
\end{equation}
In the formula \eqref{mellinQQQQ} we used the fact that
\begin{equation}
E \left\lbrace  \tilde{G}_{-\gamma}(t) \right\rbrace^{\eta -1} = \Gamma\left(\frac{1-\eta}{\gamma} + \mu \right) \frac{t^{\eta-1}}{\Gamma(\mu)}, \quad \Re \{ \eta \} < \gamma \mu +1
\end{equation}
is the Mellin transform of the density of an inverse generalized Gamma process $\tilde{G}_{-\gamma}(t)$, $t>0$ as already mentioned in the formula \eqref{mellinGgeninv}.\\
The process $\tilde{G}^1_\gamma(\tilde{G}^2_{-\gamma}(t))$, $t>0$ in the special case $\mu=\frac{1}{2}$ and $\gamma=2$ becomes a folded Cauchy process $|C(t)|$, $t>0$ and its density law has a Mellin transform given by the formula \eqref{mellinaaw} which reads
\begin{equation}
\mathcal{M}\left\lbrace q(\cdot \,,t;\frac{1}{2},2) \right\rbrace (\eta) = E \big| C(t) \big|^{\eta -1} = \frac{t^{\eta -1}}{\sin \eta \frac{\pi}{2}}, \quad 0< \Re\{ \eta \} < 2.
\label{mellinMCauchy}
\end{equation}
Hence, we can write down
\begin{equation}
|C(t)| \stackrel{i.d.}{=} \tilde{G}^1_{2,\frac{1}{2}} (\tilde{G}^2_{-2,\frac{1}{2}}(t)), \quad t>0
\label{subC}
\end{equation}
and, in view of the formula \eqref{productGtilde}, we have
\begin{equation}
|C(t)|\stackrel{i.d.}{=} \tilde{G}^1_{2,\frac{1}{2}}(t^\frac{1}{2}) \tilde{G}^2_{-2,\frac{1}{2}}(t^\frac{1}{2}) \stackrel{i.d.}{=} \frac{\tilde{G}^1_{2,\frac{1}{2}}(t^\frac{1}{2})}{\tilde{G}^2_{2,\frac{1}{2}}\left(\frac{1}{t^{1/2}} \right)}  \stackrel{i.d.}{=} \frac{1}{\tilde{G}^2_{2,\frac{1}{2}}\left(\frac{1}{t^{1/2}} \right)\tilde{G}^1_{-2,\frac{1}{2}}\left(\frac{1}{t^{1/2}} \right)}, \quad t>0.
\label{chainqq}
\end{equation}
The last equalities in distribution are due to the fact that $[\tilde{G}_\gamma(t)]^{-1} \stackrel{i.d.}{=} \tilde{G}_{-\gamma}(\frac{1}{t})$.\\
Furthermore, from the chain of equalities \eqref{chainqq}, we can state
\begin{equation}
|C(t)| \stackrel{i.d.}{=} \frac{1}{\big| C\left( \frac{1}{t} \right) \big|}, \quad t>0.
\label{iiio}
\end{equation}
Also, note that the relation \eqref{iiio} becomes 
\begin{equation}
C(t) \stackrel{i.d.}{=} \frac{1}{C\left( \frac{1}{t} \right)}, \quad t>0
\end{equation}
for the symmetry of the density law of the Cauchy process.
\begin{os}
\normalfont
From \eqref{mellinaaw} formula \eqref{mellinMCauchy} follows immediately. We derive the formula \eqref{mellinMCauchy} in a different way. It suffices to evaluate the integral
\begin{align}
\int_0^\infty x^{\eta -1} \frac{t}{\pi (t^2 + x^2)} dx = & \frac{t^{\eta -1}}{\pi} \int_0^\infty x^{\eta -1} \int_0^\infty e^{-y(1+x^2)} dy\, dx\\
= & \frac{t^{\eta -1}}{2 \pi} \Gamma\left( \frac{\eta}{2} \right)  \int_0^\infty  e^{-y} y^{-\frac{\eta}{2}} dy \nonumber \\
= & \frac{t^{\eta -1}}{2 \pi} \Gamma\left( \frac{\eta}{2} \right) \Gamma\left(1- \frac{\eta}{2} \right) \nonumber \\
= & \frac{t^{\eta -1}}{2 \sin \left( \eta \frac{\pi}{2} \right)}, \quad 2>\eta >0. \nonumber
\end{align}
Moreover, 
\begin{equation}
\int_\mathbb{R} |x|^{\eta-1} \frac{t}{\pi (t^2 + x^2)} dx = 2\int_{0}^{\infty} x^{\eta-1} \frac{t}{\pi (t^2 + x^2)} dx= \frac{t^{\eta -1}}{\sin \left( \eta \frac{\pi}{2} \right)}.
\end{equation}
\end{os}

We consider now the process $\tilde{G}^1_{-\gamma}(\tilde{G}^2_{\gamma}(t))$, $t>0$. It has distribution given by
\begin{align}
q(x,t;\mu, \gamma) = & \gamma^2 \int_{0}^{\infty} \left( \frac{x}{s} \right)^{-\gamma \mu -1} \frac{e^{-\left( \frac{x}{s} \right)^{-\gamma}}}{s \Gamma(\mu)} \left( \frac{s}{t} \right)^{\gamma \mu -1} \frac{e^{-\left( \frac{s}{t} \right)^\gamma}}{t \Gamma(\mu)} ds \nonumber \\
= & \frac{\gamma^2}{x} \frac{1}{(xt)^{\gamma \mu}} \frac{1}{\Gamma^2(\mu)} \int_0^\infty s^{2\gamma \mu -1} e^{-s^\gamma \left( \frac{1}{x^\gamma} + \frac{1}{t^\gamma} \right)} ds \nonumber \\
= & \gamma \frac{(xt)^{\gamma \mu}}{x (x^\gamma + t^\gamma)^{2\mu}} \frac{\Gamma(2\mu)}{\Gamma^2(\mu)} \nonumber
\end{align}
which coincides with \eqref{aaw}.\\
We can obtain the same result by considering the Mellin transform \eqref{mellinaaw} and the following equality in distribution
\begin{equation}
\tilde{G}^1_{-\gamma}(\tilde{G}^2_{\gamma}(t)) \stackrel{i.d.}{=} \tilde{G}^1_{-\gamma}(t^\frac{1}{2}) \tilde{G}^2_{\gamma}(t^\frac{1}{2}), \quad t>0, \; \gamma >0.
\label{productGtilde2}
\end{equation}
Indeed, from the transform \eqref{mellinGgeninv} it appears that 
\begin{equation}
E\left\lbrace \tilde{G}^1_{-\gamma}(\tilde{G}^2_{\gamma}(t)) \right\rbrace^{\eta-1} = \frac{\Gamma\left(\frac{1-\eta}{\gamma} + \mu \right)}{\Gamma(\mu)} E\left\lbrace \tilde{G}^2_\gamma(t) \right\rbrace^{\eta-1} 
\end{equation}
and we obtain back \eqref{mellinaaw}.\\
Hence, given the formulas \eqref{productGtilde2} and \eqref{productGtilde} we can conclude that
\begin{equation}
\tilde{G}^1_{-\gamma}(\tilde{G}^2_{\gamma}(t)) \stackrel{i.d.}{=} \tilde{G}^1_{\gamma}(\tilde{G}^2_{-\gamma}(t)), \quad t>0, \; \gamma >0.
\end{equation}
Furthermore, from the formula \eqref{subC}, we have
\begin{equation}
\tilde{G}^1_{-\gamma}(\tilde{G}^2_{\gamma}(t)) \stackrel{i.d.}{=} \tilde{G}^1_{-\gamma}(t^\frac{1}{2}) \tilde{G}^2_{\gamma}(t^\frac{1}{2}) \stackrel{i.d.}{=}|C(t)|, \quad t>0.
\label{CauchyRep2}
\end{equation}
for $\gamma=2$ and $\mu=\frac{1}{2}$.

The composition of two generalized Gamma processes $\tilde{G}^1_{\gamma, \mu_1}(t)$, $t>0$ and $\tilde{G}^2_{\gamma, \mu_2}(t)$, $t>0$ where the process $\tilde{G}^j_{\gamma_j, \mu_j}(t)$, $t>0$, for $j=1,2$, possesses distribution given by $\tilde{Q}(x; t, \gamma_j, \mu_j)$ are now examined.\\
The distribution of the process $\tilde{G}^1_{\gamma, \mu_1}(\tilde{G}^2_{-\gamma, \mu_2}(t))$, $t>0$, $\gamma>0$ can be easily written in an explicit form as
\begin{equation}
q(x,t)=\gamma \frac{x^{\gamma \mu_1 -1} t^{\gamma \mu_2}}{(x^\gamma + t^\gamma)^{\mu_1+\mu_2}} \frac{\Gamma(\mu_1 + \mu_2)}{\Gamma(\mu_1) \Gamma(\mu_2)}, \quad t>0, \; x >0, \; \mu_1, \mu_2>0.
\label{tdist0}
\end{equation}
For $\gamma=2$, $\mu_1=\frac{1}{2}$, $\mu_2=\frac{\nu}{2}$ and $\nu >0$ the distribution  \eqref{tdist0} becomes
\begin{equation}
f(x,t)=2 \frac{t^\nu}{(x^2 + t^2)^\frac{\nu +1}{2}} \frac{\Gamma\left( \frac{\nu+1}{2} \right)}{ \sqrt{\pi} \Gamma\left( \frac{\nu}{2} \right)}, \quad x \geq 0, \; t>0.
\label{tdist1}
\end{equation}
We point out that $f(x,\nu^\frac{1}{2})$ is the Student's t-distribution on the positive real line. \\
From the formulae \eqref{melGGen} and \eqref{mellinGgeninv} we still have
\begin{equation}
\tilde{G}^1_{\gamma, \frac{1}{2}}(\tilde{G}^2_{-\gamma, \frac{\nu}{2}}(t)) \stackrel{i.d.}{=} \tilde{G}^1_{\gamma, \frac{1}{2}}(t^\frac{1}{2}) \tilde{G}^2_{-\gamma, \frac{\nu}{2}}(t^\frac{1}{2}), \quad t>0.
\end{equation}
Indeed, the Mellin transform of the process in question reads
\begin{align}
E\left\lbrace \tilde{G}^1_{\gamma, \frac{1}{2}}(\tilde{G}^2_{-\gamma, \frac{\nu}{2}}(t)) \right\rbrace^{\eta -1} = & \frac{\Gamma\left( \frac{\eta -1}{\gamma} + \frac{1}{2} \right)}{\pi} E\left\lbrace \tilde{G}^2_{-\gamma, \frac{\nu}{2}}(t) \right\rbrace^{\eta -1} \nonumber \\
= & E\left\lbrace \tilde{G}^1_{\gamma, \frac{1}{2}}(t^\frac{1}{2}) \right\rbrace^{\eta -1} E\left\lbrace \tilde{G}^2_{-\gamma, \frac{\nu}{2}}(t^\frac{1}{2}) \right\rbrace^{\eta -1} \nonumber \\
= & E\left\lbrace \tilde{G}^1_{\gamma, \frac{1}{2}}(t^\frac{1}{2}) \, \tilde{G}^2_{-\gamma, \frac{\nu}{2}}(t^\frac{1}{2}) \right\rbrace^{\eta -1}.  \nonumber
\end{align}

\begin{os}
\normalfont 
Let us consider the process $\tilde{G}^1_{\gamma_1}(\tilde{G}^2_{\gamma_2}(t))$, $t>0$ where $\tilde{G}^1_{\gamma_1}(t)$, $t>0$ and  $\tilde{G}^2_{\gamma_2}(t)$, $t>0$ are two independent generalized Gamma processes. we have that
\begin{align*}
\left\lbrace \tilde{G}^1_{\gamma_1}(\tilde{G}^2_{\gamma_2}(t)) \, \tilde{G}^1_{\gamma_1}(\tilde{G}^2_{\gamma_2}(s)) \right\rbrace \stackrel{i.d.}{=} & \left\lbrace \tilde{G}^1_{\gamma_1}(t^\frac{1}{2}) \, \tilde{G}^2_{\gamma_2}(t^\frac{1}{2}) \, \tilde{G}^1_{\gamma_1}(s^\frac{1}{2}) \, \tilde{G}^2_{\gamma_2}(s^\frac{1}{2}) \right\rbrace
\end{align*}
For the right-hand side we are able to give the following expression for the covariance function
\begin{align*}
& E\left\lbrace \tilde{G}^1_{\gamma_1}(t^\frac{1}{2}) \, \tilde{G}^2_{\gamma_2}(t^\frac{1}{2}) \, \tilde{G}^1_{\gamma_1}(s^\frac{1}{2}) \, \tilde{G}^2_{\gamma_2}(s^\frac{1}{2}) \right\rbrace \\
&= E\left\lbrace \tilde{G}^1_{\gamma_1}(t^\frac{1}{2}) \, \tilde{G}^1_{\gamma_1}(s^\frac{1}{2}) \right\rbrace E \left\lbrace   \tilde{G}^2_{\gamma_2}(t^\frac{1}{2}) \, \tilde{G}^2_{\gamma_2}(s^\frac{1}{2}) \right\rbrace
\end{align*}
for $\gamma_1 \neq 0, \; \gamma_2 \neq 0$. 
\end{os}

\subsection{Compositions involving $n$ independent Generalized Gamma processes}

We study, in a more general setting, the composition of generalized Gamma processes. We can state the following result
\begin{te}
For $n$ independent generalized Gamma processes $\tilde{G}^j_{\gamma_j}(t)$, $t>0$, $\gamma_j \neq 0$, $j=1,2,\ldots , n$ we have that
\begin{equation}
\tilde{G}_{\gamma_1}^1(\tilde{G}_{\gamma_2}^2(\ldots \tilde{G}_{\gamma_n}^n(t) \ldots)) \stackrel{i.d.}{=} \prod_{j=1}^n \tilde{G}_{\gamma_j}^j(t^\frac{1}{n}), \quad t>0,\; |\gamma_j |>0 \, \forall j
\label{P4P}
\end{equation}
\label{teGtildeProd}
\end{te}
\paragraph{Proof:}
First of all, we write down the Mellin transform of the function $Q_j= \pm Q(x;t,\mu, \pm \gamma_j)$, $\mu>0$, $\gamma_j >0$, $j=1,2,\ldots , n$ which reads
\begin{equation}
\mathcal{M}\{ Q_j \}(\eta) = \frac{\phi_j}{\Gamma(\mu)} t^{\eta -1}
\label{wNi}
\end{equation}
where
\begin{equation}
\phi_j= \left\lbrace \begin{array}{ll} \Gamma\left( \frac{\eta -1}{\gamma_j} + \mu \right) & \textrm{if } Q_j=+Q(x;t,\mu, + \gamma_j) \\ \Gamma\left( \frac{1- \eta }{\gamma_j} + \mu \right) & \textrm{if } Q_j=-Q(x;t,\mu, -\gamma_j)  \end{array} \right . \quad j=1,2,\ldots , n.
\end{equation}
We have to highlight the presence of the term $t^{\eta -1}$ in the formula \eqref{wNi} . Such a term is closely related to the Mellin transform of the distribution of iterated processes. It turns out to be useful in the Mellin convolution machinery as well as we show below. The compound process $\tilde{G}_{\gamma_1}^1(\tilde{G}_{\gamma_2}^2(\ldots \tilde{G}_{\gamma_n}^n(t) \ldots))$ has a density law which is a composition of the density laws $Q_j$, $j=1,2,\ldots , n$. Such a composition can be written as
\begin{equation}
\circledcirc_{j=1}^n Q_j = \underbrace{\int_0^\infty \ldots \int_0^\infty}_{n-1 \textrm{ times}} Q_1(x,s_1) Q_2(s_1,s_2) ds_1 \ldots  Q_{n}(s_{n-1}, t) ds_{n-1}
\end{equation}
and thus, we have
\begin{align*}
\mathcal{M}\left\lbrace \circledcirc_{j=1}^n Q_j \right\rbrace (\eta) =  & \underbrace{\int_0^\infty \ldots \int_0^\infty}_{n-1 \textrm{ times}} \frac{\phi_1}{\Gamma(\mu)} s_1^{\eta -1} Q_2(s_1,s_2) ds_1 \ldots  Q_{n}(s_{n-1}, t) ds_{n-1} \\
= & \underbrace{\int_0^\infty \ldots \int_0^\infty}_{n-2 \textrm{ times}} \frac{\phi_1}{\Gamma(\mu)} \frac{\phi_2}{\Gamma(\mu)} s_2^{\eta -1} Q_3(s_2,s_3) ds_2 \ldots  Q_{n}(s_{n-1}, t) ds_{n-1} \\
= & \prod_{j=1}^n \left( \frac{\phi_j}{\Gamma(\mu)} t^\frac{\eta -1}{n} \right), \quad \gamma_j \neq 0 , \; \forall j.
\end{align*}
We can easily ascertain that
\begin{equation*}
E\left\lbrace \prod_{j=1}^n \tilde{G}_{\gamma_j}^j(t^\frac{1}{n}) \right\rbrace^{\eta -1} = \prod_{j=1}^n E\left\lbrace  \tilde{G}_{\gamma_j}^j(t^\frac{1}{n}) \right\rbrace^{\eta -1}= \prod_{j=1}^n \left( \frac{\phi_j}{\Gamma(\mu)} t^\frac{\eta -1}{n} \right).
\end{equation*}
Hence, 
\begin{equation*}
E \left\lbrace \tilde{G}_{\gamma_1}^1(\tilde{G}_{\gamma_2}^2(\ldots \tilde{G}_{\gamma_n}^n(t) \ldots)) \right\rbrace^{\eta -1} = E \left\lbrace  \prod_{j=1}^n \tilde{G}_{\gamma_j}^j(t^\frac{1}{n}) \right\rbrace^{\eta -1} 
\end{equation*}
and this concludes the proof. $\blacksquare$\\

The previous Theorem is a consequence of the scaling property of the independent process $\tilde{G}^i_{\gamma_i}$, $1 \leq i \leq n$. Since $\tilde{G}_\gamma(t) \stackrel{i.d.}{=}t \, \tilde{G}_\gamma(1)$ we have, tanks to the independence of the considered processes,
\begin{equation*}
\tilde{G}_{\gamma_1}^1(\tilde{G}_{\gamma_2}^2(\ldots \tilde{G}_{\gamma_n}^n(t) \ldots)) \stackrel{i.d.}{=} t \, \prod_{j=1}^n \tilde{G}_{\gamma_j}^j(1) 
\end{equation*} 
and the claimed result follows immediately.\\

The last result is related to the possibility of writing the distribution of the process $\tilde{G}_\gamma(t)$, $t>0$ (and its compositions) in terms of Fox's H-functions as shown below.\\ 
For the process $\tilde{G}_\gamma(t)$, $t>0$ with density law $Q(x;t^\gamma, \mu, \gamma)$, we have 
\begin{equation}
q(x,t) = \frac{x^{\gamma \mu -1} e^{-\frac{x^\gamma}{t^\gamma}}}{t^{\gamma \mu}\Gamma(\mu)} = \frac{(xt)^{-1}}{\Gamma(\mu)} H^{1,0}_{0,1}\left[\frac{x}{t} \Bigg| \begin{array}{c} -\\ (\mu \frac{1}{\gamma}) \end{array} \right]
\label{rap1}
\end{equation}
and, for the compound process $\tilde{G}^1_{\gamma_1}(\tilde{G}^2_{\gamma_2}(t))$, $t>0$ the distribution can be written as follows
\begin{equation}
q(x,t)=\frac{(xt)^{-1}}{\Gamma^2(\mu)} H^{2,0}_{0,2} \left[ \frac{x}{t} \Bigg| \begin{array}{cc} -;& -;\\(\mu, \frac{1}{\gamma_1});& (\mu, \frac{1}{\gamma_2}) \end{array} \right].
\label{rap2}
\end{equation}
The distribution \eqref{Gtilde} is a special case of the formula \eqref{rap2} with $\gamma_1=\gamma_2=\gamma$.\\
Finally, we give the distribution of the $n-$times iterated process \eqref{P4P} that writes
\begin{equation}
q(x,t)=\frac{(xt)^{-1}}{\Gamma^n(\mu)} H^{n,0}_{0,n} \left[\frac{x}{t} \Bigg| \begin{array}{cccc} -;& -;& \ldots ;& -\\(\mu, \frac{1}{\gamma_1});& (\mu, \frac{1}{\gamma_2});& \ldots ;& (\mu, \frac{1}{\gamma_n}) \end{array} \right].
\label{rap3}
\end{equation}
The representations \eqref{rap1} and \eqref{rap2} follow from \eqref{rap3} whereas, the last one can be verified by considering the following Mellin transform
\begin{align}
E\left\lbrace \tilde{G}^1_{\gamma_1}(\tilde{G}^2_{\gamma_2}(\ldots \right . & \left . \tilde{G}^n_{\gamma_n}(t) \ldots)) \right\rbrace^{\eta -1} \\
= & \frac{t^{\eta -1}}{\Gamma^n(\mu)} \prod_{j=1}^{n} \Gamma\left( \frac{\eta}{\gamma_j} + \mu -\frac{1}{\gamma_j} \right). \nonumber \\
= & \frac{t^{\eta-1}}{\Gamma^n(\mu)} \mathcal{M} \left\lbrace H^{n,0}_{0,n}\left[ x \Bigg| \begin{array}{cccc} -\\(\mu_j - \frac{1}{\gamma_j}, \frac{1}{\gamma_j})_{j=1,2,\ldots , n} \end{array} \right] \right\rbrace (\eta)\nonumber \\
= & [ \textrm{by } \eqref{propM1} \textrm{ and } \eqref{propH2} ]  \nonumber \\
= & \frac{t^{-1}}{\Gamma^n(\mu)} \mathcal{M} \left\lbrace \frac{1}{x} H^{n,0}_{0,n}\left[ \frac{x}{t} \Bigg| \begin{array}{cccc} -\\(\mu, \frac{1}{\gamma_j})_{j=1,2,\ldots , n} \end{array} \right] \right\rbrace (\eta). \nonumber
\end{align}

The study of the composition of Cauchy processes brings to the following interesting result
\begin{te}
For $n$ independent Cauchy processes $C(t)$, $t>0$ the following equality in distribution holds
\begin{equation}
|C^1(|C^2( \ldots |C^n(t)| \ldots)|)| \stackrel{i.d.}{=}\prod_{j=1}^n \big| C^j \left( t^\frac{1}{n} \right) \big|, \quad t>0.
\end{equation}
\end{te}
\paragraph{First proof:} It follows from \eqref{mellinMCauchy}. We evaluate, making use of the formula \eqref{mellinMCauchy}, the Mellin transform of the distribution of the folded compound process 
\[ |C^1(|C^2( \ldots |C^n(t)| \ldots)|)|, \quad  t>0 \] 
as follows
\begin{align}
E |C^1(|C^2( \ldots |C^n(t)| \ldots)|)|^{\eta-1} = & \frac{1}{\sin \left( \eta \frac{\pi}{2} \right)} E |C^2(|C^3( \ldots |C^n(t)| \ldots)|)|^{\eta-1} \label{CCCprof1} \\
= & \left\lbrace  \frac{1}{\sin \left( \eta \frac{\pi}{2} \right)} \right\rbrace^2  E |C^3(|C^4( \ldots |C^n(t)| \ldots)|)|^{\eta-1} \nonumber \\
= & \left\lbrace  \frac{1}{\sin \left( \eta \frac{\pi}{2} \right)} \right\rbrace^n t^{\eta -1} \nonumber
\end{align}
In addition, for the independence of the processes $C^j(t)$, $t>0$, $j=1,2, \ldots , n$, we have
\begin{equation*}
E\Bigg| \prod_{j=1}^n C^j \left( t^\frac{1}{n} \right) \Bigg|^{\eta-1} = \prod_{j=1}^n E\big| C^j \left( t^\frac{1}{n} \right) \big|^{\eta-1} = \left\lbrace  \frac{t^\frac{\eta -1}{n}}{\sin \left( \eta \frac{\pi}{2} \right)} \right\rbrace^n .
\end{equation*}
This concludes the first proof. 
\paragraph{Second proof:} 
For the sake of simplicity we consider only one iteration. The general case is a very trivial extension.
Making use of the representation \eqref{CauchyRep2} of the Cauchy process the following equality holds
\begin{equation}
|C^1(|C^2(t)|)| \stackrel{i.d.}{=} \tilde{G}^1_{-\gamma}(\tilde{G}^2_{\gamma}(\tilde{G}^3_{-\gamma}(\tilde{G}^4_{\gamma}(t))))
\label{Wultima}
\end{equation}
for $\gamma=2$ and $\mu=\frac{1}{2}$.\\
By applying the Theorem \ref{teGtildeProd}, we can write the formula \eqref{Wultima} as follows
\begin{equation}
\tilde{G}^1_{-\gamma}(\tilde{G}^2_{\gamma}(\tilde{G}^3_{-\gamma}(\tilde{G}^4_{\gamma}(t)))) \stackrel{i.d.}{=} \tilde{G}^1_{-\gamma}(t^\frac{1}{4}) \, \tilde{G}^2_{\gamma}(t^\frac{1}{4}) \, \tilde{G}^3_{-\gamma}(t^\frac{1}{4}) \, \tilde{G}^4_{\gamma}(t^\frac{1}{4}) , \quad t>0
\end{equation}
and thus (by \eqref{CauchyRep2})
\begin{align}
\tilde{G}^1_{-\gamma}(\tilde{G}^2_{\gamma}(\tilde{G}^3_{-\gamma}(\tilde{G}^4_{\gamma}(t)))) \stackrel{i.d.}{=} & \tilde{G}^1_{-\gamma}(\tilde{G}^2_{\gamma}(t^\frac{1}{2})) \, \tilde{G}^3_{-\gamma}(\tilde{G}^4_{\gamma}(t^\frac{1}{2})), \quad t>0 \\ 
\stackrel{i.d.}{=} & |C^1(t^\frac{1}{2})| \, |C^2(t^\frac{1}{2})|, \quad t>0. \nonumber
\end{align}
By combining such results we obtain
\begin{equation}
|C^1(|C^2(t)|)| \stackrel{i.d.}{=} |C^1(t^\frac{1}{2})| \, |C^2(t^\frac{1}{2})|
\end{equation}
and the claimed result is obtained. $\blacksquare$\\

The relation \eqref{subC} is equivalent, in some sense, to the well-known subordination law given by
\begin{equation}
 C(t) \stackrel{i.d.}{=} B(S_{\frac{1}{2}}(t)), \quad t>0
\label{equivC}
\end{equation}
where $B(t)$, $t>0$ is a standard Brownian motion and $S_{\frac{1}{2}}(t)$, $t>0$ is a $1/2-$stable law (the first passage time). In terms of the generalized Gamma process $G_\gamma(t)$, $t>0$,  we can rewrite the equivalence \eqref{equivC} as follows
\begin{equation}
| C(t) | \stackrel{i.d.}{=} G_{2,\frac{1}{2}}\left( G_{-1, \frac{1}{2}}\left(\frac{1}{t^2}\right)\right), \quad t>0.
\label{mQQQQ}
\end{equation}

\begin{os}
\normalfont
We point out that the process $G^1_\gamma(G^2_\gamma(t))$, $t>0$ is very different from $\tilde{G}^1_\gamma (\tilde{G}^2_\gamma(t))$, $t>0$. Indeed the Mellin transform of the density of $G^1_\gamma(G^2_\gamma(t))$, $t>0$ reads
\begin{equation}
E\left\lbrace G^1_\gamma(G^2_\gamma(t)) \right\rbrace^{\eta -1} = \frac{\Gamma(\frac{\eta -1}{\gamma} + \mu)}{\Gamma(\mu)} \frac{\Gamma(\frac{\eta -1}{\gamma^2} + \mu)}{\Gamma(\mu)} t^\frac{\eta-1}{\gamma^2}
\end{equation}
which is different from the formula \eqref{hhg}. Such processes have the same distribution if and only if $\gamma=1$.
\end{os}

The Cauchy process is a self-similar process in the sense that
\begin{equation*}
P\left\lbrace C(at)<w \right\rbrace = \int_{-\infty}^{w} \frac{at\, dx}{\pi (a^2t^2 + x^2)} = \int_{-\infty}^{\frac{w}{a}} \frac{t\, dz}{\pi (t^2 + z^2)} = P\left\lbrace C(t)<\frac{w}{a} \right\rbrace .
\end{equation*}
We can also use the fact that
\begin{equation}
E |aC(t)|^{\eta -1} = a^{\eta -1} E|C(t)|^{\eta-1} = \frac{(at)^{\eta -1}}{\sin \left( \eta \frac{\pi}{2} \right)} = E|C(at)|^{\eta -1}, \quad a>0.
\end{equation}
When we compose Cauchy processes, the composition maintains the self-similarity. In particular we have
\begin{align}
a|C^1(|C^2( \ldots |C^n(t)| \ldots)|)| \stackrel{i.d.}{=} & |C^1(|C^2( \ldots |C^n(at)| \ldots)|)|, \quad t>0, \; a>0 \label{m234} \\
\stackrel{i.d.}{=} & \prod_{j=1}^n \big| C^j \left( (at)^\frac{1}{n} \right) \big|, \quad t>0, \; a>0. \nonumber
\end{align}
It is enough have a look at the Mellin transform
\begin{align*}
E \big| a|C^1(|C^2( \ldots |C^n(t)| \ldots)|)| \big|^{\eta -1} = & a^{\eta -1} E|C^1(|C^2( \ldots |C^n(t)| \ldots)|)|^{\eta -1} \\
= &  [ \textrm{by } \eqref{CCCprof1} ] = \left\lbrace \frac{(at)^\frac{\eta -1}{n}}{\sin \left( \eta \frac{\pi}{2} \right)} \right\rbrace^n  \\ 
= & E \Bigg| \prod_{j=1}^n  C^j \left( (at)^\frac{1}{n} \right) \Bigg|^{\eta -1}, \quad a>0
\end{align*}
and thus, the Mellin transform of both members in the formula \eqref{m234} coincides.\\
Furthermore, for the symmetry of the Cauchy random variable,  we can write down 
\begin{equation}
a C^1(|C^2( \ldots |C^n(t)| \ldots)|) \stackrel{i.d.}{=}\prod_{j=1}^n C^j \left( (at)^\frac{1}{n} \right) , \quad t>0, \; a>0.
\end{equation}

We give now some other particular case related to the compositions of the independent generalized Gamma process $G_\gamma(t)$, $t>0$ which has density law $Q(x;t,\mu,\gamma)$, $\gamma$, $\mu>0$. 

\paragraph{The process $G^1_\gamma(G^2_1(t))$.}

The composition of a generalized Gamma process $G^1_\gamma(t)$, $t>0$ and a Gamma process $G^2_1(t)$, $t>0$ has a density law which can be written in an explicit form as
\begin{align}
q(x,t;\mu,\gamma) = & \gamma \int_0^\infty \frac{x^{\mu \gamma -1} e^{-\frac{x^\gamma}{s}}}{s^\mu \Gamma(\mu)} \frac{s^{\mu -1} e^{-\frac{s}{t}}}{t^\mu \Gamma(\mu)} ds \label{densityUNOK} \\
= & 2\frac{\gamma x^{\mu \gamma -1}}{t^\mu \Gamma^2(\mu)} K_0\left( 2 \sqrt{ \frac{x^\gamma}{t} } \right), \quad x >0 , \; t>0, \; \mu,\gamma >0. \nonumber
\end{align}
Its Mellin transform writes
\begin{equation}
E\left\lbrace G^1_\gamma(G^2_1(t)) \right\rbrace^{\eta -1} = \frac{\Gamma^2\left( \frac{\eta -1}{\gamma} +\mu \right) }{\Gamma^2(\mu)} t^\frac{\eta-1}{\gamma}, \quad \Re \{\eta \} > 1-\gamma \mu . 
\end{equation}
Such a transform can be rewritten making use of the Mellin convolution formula \eqref{conv} as
\begin{equation}
E\left\lbrace G^1_\gamma(G^2_1(t)) \right\rbrace^{\eta -1} = E \left\lbrace G^1_\gamma(t^\frac{1}{2})\right\rbrace^{\eta-1} E\left\lbrace G^2_\gamma (t^\frac{1}{2}) \right\rbrace^{\eta -1} 
\end{equation}
and thus, we can conclude that
\begin{equation}
G^1_\gamma(G^2_1(t)) \stackrel{i.d.}{=} G^1_\gamma(t^\frac{1}{2}) G^2_\gamma(t^\frac{1}{2}),  \quad t>0.
\label{GdiG}
\end{equation}
We also give the distribution of the process $G^1_{\gamma, \mu_1}(G^2_{1, \mu_2}(t))$, $t>0$ where $\mu_1 \neq \mu_2$. It can be written as follows
\begin{equation}
\label{specialGG}
q(x,t;\gamma, \mu_1, \mu_2)=\frac{2\gamma}{x\Gamma(\mu_1)\Gamma(\mu_2)} \left(\frac{x^\gamma}{t} \right)^\frac{\mu_1+\mu_2}{2} K_{\mu_2-\mu_1} \left( 2 \sqrt{\frac{x^\gamma}{t}} \right), \quad x \in D, \; t>0.
\end{equation}
From the formula \eqref{K0in0} we have
\begin{equation*}
\lim_{x \to 0^+} q(x,t;\gamma, \mu_1, \mu_2) \approx \frac{\gamma x^{\gamma \mu_2 -1}}{\Gamma(\mu_1) \Gamma(\mu_2) t^{\mu_2}} \Gamma(\mu_2 - \mu_1)
\end{equation*}
and $D=[0, \infty)$ if $\gamma \mu_2 >1$ whereas $D=(0, \infty)$ if $\gamma \mu_2 <1$.\\
Furthermore, due to the symmetry of the modified Bessel function ($K_{-\nu}=K_\nu$) we have the same result in either case $\mu_2>\mu_1$ or $\mu_1>\mu_2$ (the set $D$ will be determined by $\gamma \mu_1 \lessgtr 1$).\\ 
The equality in distribution $G^1_{\gamma, \mu_1}(G^2_{1, \mu_2}(t)) \stackrel{i.d.}{=}G^1_{\gamma, \mu_1}(t^{1/2}) \, G^2_{\gamma, \mu_2}(t^{1/2})$, $t>0$ holds as well.

\begin{os}
\normalfont
We can easily ascertain that
\begin{equation}
G^1_{2,\frac{1}{2}}(G^2_{1,\frac{1}{2}}(2t)) \stackrel{i.d.}{=} | B^1(t^{\frac{1}{2}}) B^2(t^{\frac{1}{2}}) | \quad t>0.
\label{bbN}
\end{equation}
Just calculating both Mellin transforms. We remind that $G_{2,\frac{1}{2}}(2t) \stackrel{i.d.}{=} |B(t)|$.\\
Furthermore, if the subordinator is an inverse Gamma process, then we have \eqref{mQQQQ}.
\end{os}

\begin{os}
\normalfont
We also note that $\tilde{G}^1_\gamma(\tilde{G}^2_\gamma(t)) \stackrel{i.d.}{=} G^1_\gamma(|G^2_\gamma(t^\gamma)|^\gamma)$, $t>0$. From the property $|G_\gamma(t)|^\gamma \stackrel{i.d.}{=} G_1(t)$, $t>0$ we have $G^1_\gamma(|G^2_\gamma(t^\gamma)|^\gamma) \stackrel{i.d.}{=} G^1_\gamma(G_1^2(t^\gamma))$, $t>0$.
\end{os}

For the sake of completeness we also give a representation in terms of Fox's functions of the density laws that do not have an explicit form.

\paragraph{The process $G^1_{\gamma_1}(G^2_{\gamma_2}(t))$.}

The density law of the process $G^1_{\gamma_1}(G^2_{\gamma_2}(t))$, $t>0$ writes
\begin{equation}
q(x,t;\mu,\gamma_1, \gamma_2)=\frac{x^{\mu\gamma_1 -1}}{t^\mu \Gamma^2(\mu)} \int_0^\infty s^{\mu\gamma_2 - \mu -1} e^{-\frac{x^{\gamma_1}}{s} - \frac{s^{\gamma_2}}{t}} ds, \quad x>0,\; t>0
\end{equation}
where we supposed $\mu_1=\mu_2=\mu>0$, $\gamma_1>0$, $\gamma_2>0$.\\
The Mellin transform of the function $q=q(x,t;\mu,\gamma_1, \gamma_2)$ reads
\begin{align}
\mathcal{M}\left\lbrace q \right\rbrace (\eta) = & \frac{\Gamma\left( \frac{\eta -1}{\gamma_1} + \mu \right) \Gamma\left( \frac{\eta -1}{\gamma_1 \gamma_2} +\mu \right)}{\Gamma^2(\mu)} t^{\frac{\eta -1}{\gamma_1 \gamma_2}}, \quad \eta > 1-\gamma_1 \mu.  \\
= & \frac{t^\frac{\eta -1}{\gamma_1 \gamma_2}}{\Gamma^2(\mu)} \mathcal{M}\left\lbrace H^{2,0}_{0,2}\left[ x \Bigg| \begin{array}{cc} -; & -\\ (\mu - \frac{1}{\gamma_1}, \frac{1}{\gamma_1}); & (\mu -\frac{1}{\gamma_1 \gamma_2}, \frac{1}{\gamma_1 \gamma_2}) \end{array} \right] \right\rbrace (\eta) \nonumber \\
= & [ \textrm{by } \eqref{mellinHfox} ]  \nonumber \\
= & \frac{t^{-\frac{1}{\gamma_1 \gamma_2}}}{\Gamma^2(\mu)} \mathcal{M}\left\lbrace H^{2,0}_{0,2}\left[ x t^{-\frac{1}{\gamma_1 \gamma_2}} \Bigg| \begin{array}{cc} -; & -\\ (\mu - \frac{1}{\gamma_1}, \frac{1}{\gamma_1}) ; & (\mu -\frac{1}{\gamma_1 \gamma_2}, \frac{1}{\gamma_1 \gamma_2}) \end{array} \right] \right\rbrace (\eta) \nonumber \\
= & [ \textrm{by } \eqref{propM1} ]  \nonumber \\
= & \frac{1}{t^{1/\gamma_1 \gamma_2}} \mathcal{M}\left\lbrace H^{2,0}_{2,2}\left[ \frac{x}{ t^{1/\gamma_1 \gamma_2}} \Bigg| \begin{array}{cc} (\mu, 0); & (\mu , 0)\\ (\mu - \frac{1}{\gamma_1}, \frac{1}{\gamma_1}) ; & (\mu -\frac{1}{\gamma_1 \gamma_2}, \frac{1}{\gamma_1 \gamma_2}) \end{array} \right] \right\rbrace (\eta) \nonumber \\
= & [ \textrm{by } \eqref{propH2} ]  \nonumber \\
= & \frac{1}{t^{1/\gamma_1 \gamma_2}} \mathcal{M}\left\lbrace \frac{1}{x} H^{2,0}_{2,2}\left[ \frac{x}{ t^{1/\gamma_1 \gamma_2}} \Bigg| \begin{array}{cc} (\mu, 0) ; & (\mu , 0)\\ (\mu , \frac{1}{\gamma_1}) ; & (\mu , \frac{1}{\gamma_1 \gamma_2}) \end{array} \right] \right\rbrace (\eta) \nonumber
\end{align}
and thus
\begin{equation}
q(x,t) = \frac{1}{x\, t^{1/\gamma_1 \gamma_2}} H^{2,0}_{2,2}\left[ \frac{x}{ t^{1/\gamma_1 \gamma_2}} \Bigg| \begin{array}{cc} (\mu, 0); & (\mu , 0)\\ (\mu , \frac{1}{\gamma_1}); & (\mu , \frac{1}{\gamma_1 \gamma_2}) \end{array} \right], \quad x>0,\; t>0,
\label{HGammaGamma}
\end{equation}
or
\begin{equation}
q(x,t) = \frac{1}{x\, t^{1/\gamma_1 \gamma_2} \Gamma^2(\mu)} H^{2,0}_{0,2}\left[ \frac{x}{ t^{1/\gamma_1 \gamma_2}} \Bigg| \begin{array}{cc} -;& - \\ (\mu , \frac{1}{\gamma_1});& (\mu , \frac{1}{\gamma_1 \gamma_2}) \end{array} \right], \quad x>0,\; t>0,
\end{equation}
with $\mu,\gamma_1, \gamma_2>0$.\\
In the special case where $\gamma_1=1$ and $\gamma_2=\gamma$ we have the process $G_1^1(G^2_\gamma(t))$, $t>0$ that is the composition of two independent Gamma processes $G_1^1(t)$, $t>0$ and $G_\gamma^2(t)$, $t>0$.\\
The density law of the process $G_1^1(G^2_\gamma(t))$, $t>0$ can be written as
\begin{equation}
q(x,t;\gamma)=\frac{\gamma x^{\mu-1}}{t^\mu \Gamma^2(\mu)} \int_0^\infty s^{\mu\gamma -\mu -1} e^{-\frac{x}{s} - \frac{s^\gamma}{t}} ds, \quad x>0,\; t>0, \; \mu,\gamma>0.
\label{densDD}
\end{equation}
Such a function has the Mellin transform
\begin{equation}
\mathcal{M}\left\lbrace q(\cdot \, ,t;\gamma) \right\rbrace (\eta) = \frac{\Gamma(\eta +\mu -1)}{\Gamma(\mu)} \frac{\Gamma\left( \frac{\eta -1}{\gamma} +\mu \right)}{\Gamma(\mu)} t^{\frac{\eta -1}{\gamma}}, \quad \Re \{ \eta \} > 1-\mu
\label{mellinDD}
\end{equation}
Hence, the function \eqref{densDD} can be expressed in terms of Fox's H-functions by performing the inverse Mellin transform of the formula \eqref{mellinDD}. We have
\begin{align}
q(x,t;\gamma) = & \frac{t^{-\frac{1}{\gamma}}}{\Gamma^2(\mu)} \frac{1}{2\pi  i} \int_{\theta - i \infty}^{\theta + i \infty} \Gamma(\eta +\mu -1) \Gamma\left( \frac{\eta -1}{\gamma} +\mu \right) t^{\frac{\eta}{\gamma}} x^{-\eta} d\eta \\
= & \frac{1}{x \, t^{1/\gamma}} H^{2,0}_{2,2} \left[ \frac{x}{t^{1/\gamma}} \bigg| \begin{array}{cc} (\mu , 0); & (\mu , 0) \\ (\mu, 1) ; & (\mu, \frac{1}{\gamma}) \end{array} \right], \quad x>0,\; t>0 \label{dfH}
\end{align}
where $\theta$ is a real number such that $\theta > 1-\mu$.\\
The function \eqref{dfH} can be obtained by the formula \eqref{HGammaGamma}.\\

By combining a generalized Gamma process $G_\gamma^1(t)$, $t>0$ and an inverse generalized Gamma process $G^2_{-\gamma}(t)$, $t>0$ we obtain the process $G^1_\gamma(G^2_{-\gamma}(t))$, $t>0$ with distribution given by
\begin{equation}
q(x,t;\mu,\gamma)=\frac{x^{\mu\gamma -1}}{t^\mu \Gamma^2(\mu)} \int_0^\infty s^{-\mu -\gamma \mu -1}e^{-\frac{x^\gamma}{s}} e^{-\frac{1}{s^{\gamma}t}} ds, \quad x >0, t>0, \; \mu,\gamma>0.
\label{densTT}
\end{equation}
The Mellin transform of the density \eqref{densTT} reads
\begin{equation}
\mathcal{M}\left\lbrace q(\cdot \,,t;\mu,\gamma) \right\rbrace (\eta) = \frac{\Gamma\left( \frac{\eta -1}{\gamma} +\mu \right) \Gamma\left( \mu - \frac{\eta -1}{\gamma^2} \right)}{\Gamma^2(\mu)} t^{-\frac{\eta -1}{\gamma^2}}, \quad 1-\mu\gamma < \eta < \mu \gamma^2 +1.
\end{equation}
For $\gamma=1$ and $\mu=1/2$, we have
\begin{equation}
\Gamma(\eta-\frac{1}{2}) \Gamma(1- [ \eta -\frac{1}{2}]) \frac{t^{1-\eta}}{\pi} = \frac{t^{1-\eta}}{\sin( \frac{\pi}{2}[\eta - \frac{1}{2}] ) }, \quad  \frac{1}{2} < \eta < \frac{3}{2}
\end{equation}
which is the Mellin transform of the distribution of the process  $|C(\frac{1}{\sqrt{t}})|^2$. This follows from the equality \eqref{mQQQQ} by noting that $[G_\gamma(t)]^\gamma \stackrel{i.d.}{=} G_1(t)$.\\
We can give a representation of the function \eqref{densTT} in terms of H-functions as follows
\begin{align}
\mathcal{M} \left\lbrace q(\cdot \,,t;\mu,\gamma) \right\rbrace = & \frac{\Gamma\left( \mu -\frac{1}{\gamma} + \frac{\eta}{\gamma}  \right) \Gamma\left( \mu + \frac{1}{\gamma^2} - \frac{\eta}{\gamma^2} \right)}{\Gamma^2(\mu)} t^{-\frac{\eta -1}{\gamma^2}} \nonumber \\
= & \frac{t^{-\frac{\eta -1}{\gamma^2}}}{\Gamma^2(\mu)} \mathcal{M} \left\lbrace H^{1,1}_{1,1} \left[ x \Bigg| \begin{array}{c} (1-\mu-\frac{1}{\gamma^2}, \frac{1}{\gamma^2}) \\ (\mu -\frac{1}{\gamma}, \frac{1}{\gamma}) \end{array} \right] \right\rbrace (\eta) \nonumber\\
= & [ \textrm{by } \eqref{propM1} \textrm{ and } \eqref{propH2}] \nonumber \\
= & \frac{t^\frac{1}{\gamma^2}}{\Gamma^2(\mu)} \mathcal{M} \left\lbrace \frac{1}{x} H^{1,1}_{1,1} \left[ x\, t^\frac{1}{\gamma^2}  \bigg| \begin{array}{c} (1-\mu, \frac{1}{\gamma^2}) \\ (\mu , \frac{1}{\gamma}) \end{array} \right] \right\rbrace (\eta)  \nonumber
\end{align}
and thus
\begin{equation}
q(x,t;\mu,\gamma) = \frac{t^\frac{1}{\gamma^2}}{\Gamma^2(\mu)} \frac{1}{x} H^{1,1}_{1,1} \left[ x\, t^\frac{1}{\gamma^2}  \bigg| \begin{array}{c} (1-\mu, \frac{1}{\gamma^2}) \\ (\mu , \frac{1}{\gamma}) \end{array} \right], \quad x>0, t>0
\end{equation}
or equivalently
\begin{equation}
q(x,t;\mu,\gamma) = \frac{t^\frac{1}{\gamma^2}}{x} H^{1,1}_{3,1} \left[ x\, t^\frac{1}{\gamma^2}  \bigg| \begin{array}{ccc} (1-\mu, \frac{1}{\gamma^2});&(\mu,0);&(\mu,0) \\ (\mu , \frac{1}{\gamma});&-;&- \end{array} \right] , \quad x>0, t>0
\end{equation}
with $\mu>0$, $\gamma>0$.\\

We dwell now on the composition 
\begin{equation}
I_G^{n}(t)= G^1_{\gamma_1} (G^2_{\gamma_2}( \ldots G^{n+1}_{\gamma_{n+1}}(t) \ldots)), \quad t>0
\label{IGn}
\end{equation}
where $G^j_{\gamma_j}(t)$, $t>0$ for $j=1,2,\ldots , n+1$ are independent Gamma processes with density $Q(x;t,\mu,\gamma_j)$. The most important tool in the study of such a process, once again, is the Mellin transform. For the density of the process $I_G^{n}(t)$, $t>0$, say $q(x,t)$, with $\mu_j=\mu>0$, $\gamma_j>0$, $j=1,2, \ldots , n$ we have
\begin{equation}
\mathcal{M}\left\lbrace q(\cdot \,,t) \right\rbrace (\eta) = \frac{\Gamma\left( \frac{\eta -1}{\gamma_1}  + \mu \right) \Gamma\left(\frac{\eta - 1}{\gamma_1 \gamma_2} + \mu \right)}{\Gamma(\mu) \Gamma(\mu)} \cdots \frac{\Gamma\left( \frac{\eta -1}{\prod_{i=1}^{n} \gamma_i} + \mu \right)}{\Gamma(\mu)} t^{\frac{\eta -1}{\prod_{i=1}^{n} \gamma_i}}.
\label{mellinGn}
\end{equation}
The transform \eqref{mellinGn} contains fundamental informations about the structure of the moments and the density of the process  $I_G^{n}(t)$, $t>0$, in particular we can state the following theorem
\begin{te}
Let us consider the process $I_G^{n-1}(t)= G^1_{\gamma_1} (G^2_{\gamma_2}( \ldots G^{n}_{\gamma_{n}}(t) \ldots))$, $t>0$ where $G^j_{\gamma_j}(t)$, $t>0$, $j=1,2,\ldots , n$ are independent generalized Gamma processes. If $\gamma_1=\gamma > 0$ and  $\gamma_2=\gamma_3= \ldots =\gamma_n=1$, then the following equality in distribution holds
\begin{equation}
G^1_{\gamma}(G^2_1(\ldots (G^n_1(t)) \ldots)) \stackrel{i.d.}{=} \prod_{j=1}^n G_{\gamma}^j(t^\frac{1}{n}) ,\quad t>0.
\label{prodGn}
\end{equation} 
\label{teP1}
\end{te}
\paragraph{Proof:}
Fixing $\gamma_2=\gamma_3= \ldots =\gamma_n=1$ in the formula \eqref{mellinGn}, we obtain the Mellin transform of the density law of the process $I_G^{n-1}(t)$, $t>0$ which reads
\begin{equation}
\mathcal{M}\left\lbrace q \right\rbrace (\eta) = \left\lbrace \frac{\Gamma\left( \frac{\eta -1}{\gamma_1}  + \mu \right)}{\Gamma(\mu)}  t^{\frac{\eta -1}{n \gamma_1}}\right\rbrace^n .
\label{meLLp}
\end{equation}
We immediately recognize the structure of the product of independent random variables.\\
The processes $G^j_{\gamma_1}(t)$, $t>0$, $j=1,2,\ldots , n$ are independent and thus
\begin{equation*}
E\left\lbrace \prod_{j=1}^n  G^j_{\gamma_1}(t) \right\rbrace^{\eta -1} = \prod_{j=1}^n E\left\lbrace G^j_{\gamma_1}(t) \right\rbrace^{\eta -1}=\prod_{j=1}^n \frac{\Gamma\left( \frac{\eta -1}{\gamma_1}  + \mu \right)}{\Gamma(\mu)}  t^{\frac{\eta -1}{\gamma_1}}.
\end{equation*}
The last formula permits us to write down
\begin{equation*}
E\left\lbrace G^1_{\gamma_1}(G^2_1(\ldots (G^n_1(t)) \ldots)) \right\rbrace^{\eta -1}  =  E\left\lbrace \prod_{j=1}^n  G^j_{\gamma_1}(t^\frac{1}{n}) \right\rbrace^{\eta -1}
\end{equation*}
and this concludes the proof. $\blacksquare$\\

In a similar way as in the previous proof, we can show that the process \eqref{prodGn} is $\gamma_1-$self-similar in the sense that
\begin{equation}
a^\frac{1}{\gamma_1} \left\lbrace G^1_{\gamma_1}(G^2_1(\ldots (G^n_1(t)) \ldots)) \right\rbrace  \stackrel{i.d.}{=} \prod_{j=1}^n  G^j_{\gamma_j}( (at)^\frac{1}{n}), \quad a>0, \; \gamma_1 >0.
\label{selfG}
\end{equation}
Indeed, we have
\begin{align*}
E\left\lbrace a^\frac{1}{\gamma_1} G^1_{\gamma_1}(G^2_1(\ldots (G^n_1(t)) \ldots)) \right\rbrace^{\eta -1} = & a^\frac{\eta -1}{\gamma_1} E \left\lbrace G^1_{\gamma_1}(G^2_1(\ldots (G^n_1(t)) \ldots)) \right\rbrace ^{\eta -1}\\
=  & a^\frac{\eta -1}{\gamma_1} \left\lbrace \frac{\Gamma\left( \frac{\eta -1}{\gamma_1}  + \mu \right)}{\Gamma(\mu)} \right\rbrace^n t^\frac{\eta -1}{\gamma_1}\\
= & [ \textrm{by } \eqref{meLLp} ]  \\
= & E\left\lbrace \prod_{j=1}^n  G^j_{\gamma_j}( (at)^\frac{1}{n}) \right\rbrace^{\eta -1}, \quad a>0.
\end{align*}
Furthermore, we observe that
\begin{equation}
G_{\gamma_1}(G^1_{1}(G^2_1(\ldots (G^n_1(t)) \ldots)) ) \stackrel{i.d.}{=} G_{\gamma_1} \left( \prod_{j=1}^n G^j_1(t^\frac{1}{n}) \right), \quad t>0
\label{selfG1}
\end{equation}
where we recall that $G_{\gamma_1}(t) \sim Q(x;t, \mu, \gamma_1)$.\\ 

For the process $G^1_{\gamma_1} (G^2_{\gamma_2}( \ldots G^{n}_{\gamma_{n}}(t) \ldots))$, $t>0$ where $G^j_{\gamma_j}(t)$, $t>0$ for $j=1,2,\ldots n$ are independent Gamma processes with distribution $Q(x;t,\mu_j,\gamma_j)$, $\gamma_j \neq 0$, $\mu_j>0$, $j=1,2,\ldots , n$ the following equality in distribution holds
\begin{equation}
\left\lbrace G^1_{\gamma_1} (G^2_{\gamma_2}( \ldots G^{n}_{\gamma_{n}}(t) \ldots)) \right\rbrace^{\gamma_1}  \stackrel{i.d.}{=} G^1_{1} (G^2_{\gamma_2}( \ldots G^{n}_{\gamma_{n}}(t) \ldots)), \quad t>0, \; \gamma_1 \neq 0.
\label{selfGcomp}
\end{equation}
It suffices to consider the process $Z(t)=G_\gamma(X(t))$, $t>0$ with density law given by
\begin{equation}
f_Z(x,t) = \int_{0}^\infty Q(x;c, \mu, \gamma) f_X(c,t)dc, \quad x>0, \; t>0
\end{equation} 
for some process $X(t)$, $t>0$ with distribution $f_X(x,t)$.\\
We have already seen that $[G_\gamma(t)]^\gamma \stackrel{i.d.}{=}G_1(t)$ for $\gamma >0$ and $[G_\gamma(t)]^{-1} \stackrel{i.d.}{=} G_{-\gamma}(t)$ for $\gamma >0$. By combining such results we obtain $[G_\gamma(t)]^\gamma \stackrel{i.d.}{=}G_1(t)$ for $\gamma \in \mathbb{R} \setminus \{ 0 \}$. The density of the process $G_1(t)$ is $Q(x;t,\mu,1)$ and then the process $[Z(t)]^\gamma$, $t>0$ has a density law which reads
\begin{equation}
f_{Z^\gamma}(x,t) = \int_{0}^\infty Q(x;c, \mu, 1) f_X(c,t)dc \quad x>0, \, t>0.
\label{eqSS}
\end{equation}
The formula \eqref{eqSS} is the density law of the process $G_1(X(t))=[Z(t)]^\gamma$, $t>0$ for any process $X(t)$, $t>0$.\\

It is also straightforward to ascertain that
\begin{align*}
\left\lbrace G_{\gamma}(G^1_{1}(G^2_1(\ldots (G^n_1(t)) \ldots)) ) \right\rbrace^\beta  \stackrel{i.d.}{=} & G_{\frac{\gamma}{\beta}} \left( \prod_{j=1}^n G^j_1(t^\frac{1}{n}) \right), \quad t>0 \\ 
\stackrel{i.d.}{=} & G_{\frac{\gamma}{\beta}}(t^\frac{1}{n+1}) \prod_{j=1}^n G^j_1(t^\frac{1}{n+1}), \quad t>0 
\end{align*}
for $\gamma$,$\beta \in \mathbb{R}\setminus \{ 0 \}$.\\

It must be noted that the shape parameter $\gamma_1$ turns out to be very important and one can ascribe such a matter to the fact that $G^1_{\gamma_1}(t)$, $t>0$ is the guiding process.
 
\begin{os}
\normalfont
We have already pointed out that $\tilde{G}^1_\gamma(\tilde{G}^2_\gamma(t)) \stackrel{i.d.}{=} G^1_\gamma(|G^2_\gamma(t^\gamma)|^\gamma) $ and thus
\begin{align*}
\tilde{G}^1_\gamma( \tilde{G}^2_\gamma(\ldots \tilde{G}^n_\gamma(t) \ldots)) \stackrel{i.d.}{=} & G^1_\gamma(|G^2_\gamma( \ldots |G^n_\gamma(t^\gamma)|^\gamma \ldots |^\gamma), \quad t>0\\
\stackrel{i.d.}{=} & G^1_\gamma(G_1^2( \ldots G^n_1(t^\gamma) \ldots )) , \quad t>0 \\
\stackrel{i.d.}{=} & G^1_\gamma \left( \prod_{j=2}^n G_1^j(t^\frac{\gamma}{n-1})  \right), \quad t>0 \\
\stackrel{i.d.}{=} & \prod_{j=1}^n G_\gamma^j(t^\frac{\gamma}{n}) \stackrel{i.d.}{=} \prod_{j=1}^n \tilde{G}_\gamma^j(t^\frac{1}{n})  \quad t>0
\end{align*}
where we used the fact that $|G_\gamma(t)|^\gamma \stackrel{i.d.}{=} G_1(t)$ and the Theorem \ref{teP1}. 
\end{os}

We give now an interesting representation in terms of H-functions of the distribution of the process $G_1(t)$, $t>0$ and its compositions. We readily have that
\begin{equation}
Q(x;t,\mu,1) = \frac{x^{\mu -1} e^{-\frac{x}{t}}}{t^\mu \Gamma(\mu)} = \frac{x^{\mu -1}}{t^\mu \Gamma(\mu)} H^{1,0}_{0,1}\left[ \frac{x}{t} \Bigg| \begin{array}{c} -\\ (0,1)  \end{array} \right]
\end{equation}
where $Q(x;t,\mu,1)$ is the distribution of the process $G_1(t)$, $t>0$.\\
We are also able to write down the distribution of the compound process $G^1_1(G^2_1(t))$, $t>0$ as
\begin{equation}
q(x,t)=2\frac{x^{\mu-1}}{t^\mu \Gamma(\mu)} K_0\left( 2 \sqrt{\frac{x}{t}} \right) = \frac{x^{\mu -1}}{t^\mu \Gamma^2(\mu)} H^{2,0}_{0,2} \left[ \frac{x}{t} \Bigg| \begin{array}{ccc} -&;&-\\(0,1)&;&(0,1) \end{array} \right].
\end{equation}
In general, for the $n-$times iterated process, we have
\begin{align}
q(x,t)= & \frac{x^{\mu -1}}{t^\mu \Gamma^n(\mu)} H^{n,0}_{0,n}\left[ \frac{x}{t} \Bigg| \begin{array}{cccc} -;&-;&\ldots ;&-\\ (0,1)_1;&(0,1)_2;& \ldots ;& (0,1)_n  \end{array} \right] \label{GMei}\\
= & \frac{x^{\mu -1}}{t^\mu \Gamma^n(\mu)} G^{n,0}_{0,n}\left[ \frac{x}{t} \Bigg| \begin{array}{cccc} -\\ \textbf{0}  \end{array} \right], \quad x>0, \; t>0 \nonumber
\end{align}
where $G^{n,0}_{0,n}(z)$ is the Meijer's G-function (for the definition, see Mathai and Saxena \cite{MS78}).\\
Furthermore,  the distribution of the composition \eqref{prodGn} can be written in terms of H-functions as follows
\begin{equation}
q(x,t)=\frac{x^{\gamma \mu -1}}{t^\mu \Gamma^n(\mu)} H^{n,0}_{0,n}\left[ \frac{x}{t^{1/\gamma}} \Bigg| \begin{array}{cccc} -;&-;&\ldots ;&-\\ (0,\frac{1}{\gamma})_1;&(0,\frac{1}{\gamma})_2;& \ldots ;& (0,\frac{1}{\gamma})_n  \end{array} \right], \quad x>0, \; t>0.
\label{jjjL}
\end{equation} 
It suffices to evaluate the Mellin transform of the function \eqref{jjjL}. First of all, consider
\begin{equation}
\mathcal{M} \left\lbrace q(\cdot ,t) \right\rbrace (\eta) = \frac{1}{t^\mu \Gamma^n(\mu)} \mathcal{M} \left\lbrace H^{n,0}_{0,n}\left[ \frac{x}{t^{1/\gamma}} \Bigg| \begin{array}{c} -;\\ (0,\frac{1}{\gamma})_{i=1,2,\ldots ,n}  \end{array} \right] \right\rbrace (\eta + \gamma \mu -1) 
\end{equation}
by the property \eqref{propM2} and,
\begin{equation}
\mathcal{M} \left\lbrace q(\cdot \,,t) \right\rbrace (\eta) = \frac{t^\frac{\eta-1}{\gamma}}{\Gamma^n(\mu)} \mathcal{M} \left\lbrace H^{n,0}_{0,n}\left[ x \Bigg| \begin{array}{c} -;\\ (0,\frac{1}{\gamma})_{i=1,2,\ldots ,n}  \end{array} \right] \right\rbrace (\eta + \gamma \mu -1) 
\label{stepMGn}
\end{equation}
by the property \eqref{propM1}. We obtain the claimed result by observing that
\begin{equation}
\mathcal{M} \left\lbrace H^{n,0}_{0,n}\left[ \cdot  \Bigg| \begin{array}{c} -;\\ (0,\frac{1}{\gamma})_{i=1,2,\ldots ,n}  \end{array} \right] \right\rbrace (\eta + \gamma \mu -1) = \prod_{i=1}^n \Gamma\left( \frac{\eta + \gamma \mu -1}{\gamma} \right)
\end{equation}
and the formula \eqref{stepMGn} coincides with the Mellin transform \eqref{mellinGn} with $\gamma_2=\gamma_3= \ldots , \gamma_n=1$ and $\gamma_1=\gamma$.\\
The representation \eqref{GMei} follows from \eqref{jjjL}.\\

The Gamma process $G_1(t)$, $t>0$ possesses a density that is an infinitely divisible law and its Laplace transform reads
\begin{equation}
\mathcal{L}\left\lbrace \frac{x^{\mu-1} e^{-\frac{x}{t}}}{t^\mu \Gamma(\mu)} \right\rbrace = \frac{1}{(1+t\lambda)^\mu}.
\end{equation}
We exploit now the infinitely divisibility of the distribution of the process $G_1(t)$, $t>0$. It is useful to remind that for a Gamma process $G_1(t)$, $t>0$ and a generalized Gamma process $G_\gamma(t)$, $t>0$ the following hold
\begin{equation}
\left[ G_1(t) \right]^{1/\gamma} \stackrel{i.d.}{=} G_\gamma(t) \quad \textrm{and} \quad \left[ G_\gamma(t) \right]^{\gamma} \stackrel{i.d.}{=} G_1(t), \quad t>0.
\label{propertiesG}
\end{equation}
The sum of $n$ independent Gamma processes still has Gamma distribution being the Gamma density an infinitely divisible law. In light of the properties \eqref{propertiesG} we can consider a sum involving $n$ independent generalized Gamma processes $G^i_\gamma(t)$, $t>0$ $i=1,2, \ldots , n$ and we are able to write down the following equality
\begin{equation}
\sum_{i=1}^n \left[ G^i_{\gamma_i, \mu} (t) \right]^{\gamma_i} \stackrel{i.d.}{=} \sum_{i=1}^{n} G^i_{1, \mu}(t) \stackrel{i.d.}{=} G_{1,n \mu}(t), \quad t>0
\label{propSomma}
\end{equation}
and $\mu>0$, $ \gamma_i \neq 0$ for $i=1,2,\ldots, n$.\\
The properties \eqref{propertiesG}  also allow us to write the following relation
\begin{equation}
\left\lbrace  \sum_{i=1}^n \left[ G^i_{\gamma_i, \mu} (t) \right]^{\gamma_i} \right\rbrace^{1/\gamma}  \stackrel{i.d.}{=} G_{\gamma,n \mu}(t), \quad t>0
\label{GenG}
\end{equation}
and thus we can generalize the Bessel process. The density law of the process \eqref{GenG} solves the p.d.e.
\begin{equation}
\frac{\partial}{\partial t} Q = \frac{1}{\gamma^2} \left\lbrace \frac{\partial}{\partial x} x^{2-\gamma} \frac{\partial}{\partial x} Q - (\gamma n \mu -1) \frac{\partial}{\partial x} x^{1-\gamma} Q \right\rbrace, \quad x >0, \; t>0, n \in \mathbb{N}
\end{equation}
as we shown in the Theorem \ref{teCharming}.\\
A straightforward check is the case where $n$ processes $G_{2,\frac{1}{2}}(t)$, $t>0$ are involved. We obtain
\begin{equation}
\sqrt{ \sum_{i=1}^n \left[ B^i(t) \right]^2 } \stackrel{i.d.}{=} G_{2, \frac{n}{2}}(2t), \quad t>0
\end{equation}
which is the Bessel process starting from zero.

\section{Compositions involving Generalized Gamma process and Brownian motion}

In this section we will introduce and study some iterated processes. Several authors studied the compositions of processes, the most popular one is probably the iterated Brownian motion which has been introduced by Burdzy in 1993 \cite{BU931} and then studied thoroughly by Burdzy \cite{BU932}, Khoshnevisian and Lewis \cite{KL96} and other authors. The study of the iterated Brownian motion has been motivated by the analysis of diffusions in cracks (see Chudnovsky and Kunin \cite{CK87}, Kunin and Gorelik \cite{KG91}). The iterated folded Brownian motion $|B^1(|B^2(t)|)|$, $t>0$, where $B^j(t)$, $t>0$, $j=1,2$ are independent Brownian motions, has the same distribution of the process $G_{2,\frac{1}{2}}^1(G^2_{2, \frac{1}{2}}(2t))$, $t>0$ and the distribution can be rewritten in terms of H-functions by looking at the formula \eqref{HGammaGamma}. The iterated folded fractional Brownian motion $|B^1_{H_1}(|B^2_{H_2}(t)|)|$, $t>0$ where $B^j_{H_j}(t)$, $t>0$ for $j=1,2,$ are independent fractional Brownian motions, possesses the same distribution of the process $|B(G_{\frac{1}{H_1},\frac{1}{2}}(2t^{2H_2}))|$, $t>0$ and can be seen as the composition $G^1_{2,\frac{1}{2}} (G^2_{\frac{1}{H_1},\frac{1}{2}}(2t^{2H_2}))$, $t>0$ involving two independent generalized Gamma processes. Such a distribution can be written either in terms of H-functions as in the formula \eqref{HGammaGamma} or in the form
\begin{equation}
u(a,b,\gamma)=\int_0^\infty s^{\gamma -1} e^{-as - bs^{-\rho}} ds, \qquad a,b,\rho>0.
\label{termo}
\end{equation}
Several authors (see for instance Mathai and Haubold \cite{MH08}) have thoroughly studied the function \eqref{termo} and its representation in terms of $H-$functions. The function \eqref{termo} is the standard thermonuclear function in the Maxwell-Boltzmann case, in the theory of nuclear reactions. It appears very important in physics and astrophysics when we consider the probability for a thermonuclear reaction to occur in the solar fusion plasma.\\
Generally speaking, the results shown in the previous section about the generalized Gamma processes and their compositions can be useful to study a number of well-known processes.\\
In the literature, there are a lot of papers devoted to the study of the variance Gamma process where a random time change is given by a Gamma process $\Gamma(t)$, $t>0$ with distribution given by the formula \eqref{gammaLevy} (see Madam, Carr and Chang \cite{MCC}, Madan and Seneta \cite{MS90}, Kozubowski, Meerschaert and Pod\`orski \cite{KMP06}). The variace Gamma process is also termed Laplace motion since the increments follow the Laplace distribution. The variance Gamma process $B(\Gamma(t))$, $t>0$ has Laplace transform 
\begin{equation}
\mathcal{L}\{ q \}(\lambda) = \frac{1}{(1+\lambda)^t}
\end{equation} 
where $q$, once again, is that in \eqref{gammaLevy}. It is well-known that the following representation is possible
\begin{equation*}
B(\Gamma(t)) \stackrel{i.d.}{=} \Gamma_1(t) + \Gamma_2(t), \quad t>0
\end{equation*}
where $\Gamma_1$, $\Gamma_2$ are two independent Gamma subordinators with different parameters.\\
Here we study the compositions involving the Gamma process $G_1(t)$, $t>0$ and in particular we begin with the compositions where it only appears as random time. For the process $B(G_1(t))$, $t>0$ where $B(t)$, $t>0$ is a standard Brownian motion and $G_1(t)$, $t>0$ is a Gamma process with distribution $Q(x;t,\mu,1)$, $\mu>0$, we give the following result
\begin{te}
The process $B(G_1(t))$, $t>0$ has a distribution given by
\begin{align}
q(x,t;\mu) = & \int_0^\infty \frac{e^{-\frac{x^2}{2s}}}{\sqrt{2\pi s}} Q(s;t,\mu,1) ds \label{densityBG1}\\
= & 2\frac{|x|^{\mu -\frac{1}{2}} }{ \pi^{\frac{1}{2}} \Gamma(\mu)} \left(\frac{2}{t} \right)^\frac{2\mu +1}{4} K_{\mu -\frac{1}{2}} \left( |x| \sqrt{\frac{2}{t}} \right), \quad x \in \mathbb{R},\; t>0, \; \mu>0 \nonumber
\end{align}
which solves the following p.d.e.
\begin{align}
\frac{\partial}{\partial t} q = & \frac{\mu}{2} \frac{\partial^2}{\partial x^2}q - \frac{1}{4} \frac{\partial^3}{\partial x^3} \left(x\,q \right) \label{operatorQQ} \\
= & \frac{2\mu -3}{4} \frac{\partial^2}{\partial x^2} q - \frac{x}{4} \frac{\partial^3}{\partial x^3} q, \quad x\in \mathbb{R}, \; t>0, \; \mu \in \left( 0, \frac{1}{2}\right) \cup \left( \frac{1}{2}, \infty \right) \nonumber
\end{align}
where $q=q(x,t;\mu)$.
\label{teQQ}
\end{te}
\paragraph{First proof:} From the property \eqref{K0in0} of the modified Bessel function we can observe that
\begin{equation*}
\lim_{x \to 0^+} q(x,t;\mu) = \frac{2^{\frac{\mu}{4}}}{ \sqrt{\pi}\Gamma(\mu)} \left(\frac{1}{t} \right)^{\frac{\mu}{4}+\frac{1}{2}}\Gamma\left( \mu -\frac{1}{2}\right) \quad t>0, \; \mu>0, \quad \left( \mu \neq \frac{1}{2} \right)
\end{equation*}
and
\begin{equation*}
\lim_{x \to 0^+} q(x,t;\frac{1}{2}) = \infty.
\end{equation*}
Consider now the functions
\begin{equation*}
p(x,s)=\frac{e^{-\frac{x^2}{2s}}}{\sqrt{2\pi s}}\quad \textrm{ and  } \quad Q(s,t)=Q(s;t,\mu,1).
\end{equation*}
We evaluate the time derivative
\begin{align*}
\frac{\partial}{\partial t} q = & [ \textrm{ by } \eqref{PDEm1} ] = \int_0^\infty p(x,s) \left\lbrace s \frac{\partial^2}{\partial s^2} - (\mu-2) \frac{\partial}{\partial s} \right\rbrace Q(s,t) ds\\
= & \int_0^\infty p(x,s) s \frac{\partial^2}{\partial s^2} Q ds - (\mu -2) \int_0^\infty p(x,s) \frac{\partial}{\partial s} Q(s,t) ds.
\end{align*}
By performing an integration by parts one obtains
\begin{align*}
\frac{\partial}{\partial t} q = &  p(x,s) s \frac{\partial}{\partial s} Q(s,t) \Bigg|_{s=0}^{s=\infty} - \int_0^\infty \frac{\partial}{\partial s} \left\lbrace p(x,s) s \right\rbrace \frac{\partial}{\partial s} Q(s,t) ds \\
- & (\mu -2) \int_0^\infty p(x,s) \frac{\partial}{\partial s} Q(s,t) ds \\
= & p(x,s) s \frac{\partial}{\partial s} Q(s,t) \Bigg|_{s=0}^{s=\infty} - \int_0^\infty s \frac{\partial}{\partial s} p(x,s) \frac{\partial}{\partial s} Q(s,t) ds \\ 
- & (\mu -1) \int_0^\infty p(x,s) \frac{\partial}{\partial s} Q(s,t) ds
\end{align*}
where the coefficients of the integral $\int_0^\infty pQ^\prime ds$ have been summed. We integrate by parts once again,
\begin{align*}
\frac{\partial}{\partial t} q = & p(x,s) s \frac{\partial}{\partial s} Q(s,t) \Bigg|_{s=0}^{s=\infty} - s Q(s,t) \frac{\partial}{\partial s}p(x,s) \Bigg|_{s=0}^{s=\infty} + \int_0^\infty \frac{\partial}{\partial s}\left(p(x,s) \right) Q(s,t) ds  \\
+ & \int_0^\infty s \frac{\partial^2}{\partial s^2} \left( p(x,s) \right) Q(s,t)ds - (\mu -1) \int_0^\infty p(x,s) \frac{\partial}{\partial s} Q(s,t) ds.
\end{align*}
Also we observe that
\begin{align*}
 (\mu-1)\int_0^\infty p(x,s) \frac{\partial}{\partial s} Q(s,t) ds = & (\mu-1)  \left\lbrace p(x,s)Q(s,t)\Bigg|_{s=0}^{s=\infty}  \right . \\
& \left .-  \int_0^\infty \frac{\partial}{\partial s}\left(p(x,s) \right) Q(s,t) ds \right\rbrace 
\end{align*}
and
\begin{align*}
p(x,s) s \frac{\partial}{\partial s} Q(s,t) \Bigg|_{s=0}^{s=\infty} = & p(x,s)\left\lbrace (\mu-1) - \frac{s}{2t} \right\rbrace Q(s,t) \Bigg|_{s=0}^{s=\infty} \\
= & (\mu -1)p(x,s) Q(s,t) \Bigg|_{s=0}^{s=\infty}
\end{align*}
because
\begin{equation*}
p(x,s) \frac{s}{2t}Q(s,t) \Bigg|_{s=0}^{s=\infty} = \frac{e^{-\frac{x^2}{2s}}}{\sqrt{2\pi s}} \frac{1}{2t} \frac{s^\mu e^{-\frac{s}{t}}}{t^\mu \Gamma(\mu)} \Bigg|_{s=0}^{s=\infty} =0, \qquad (\mu >0).
\end{equation*}
Moreover,
\begin{equation*}
 - s Q(s,t) \frac{\partial}{\partial s}p(x,s) \Bigg|_{s=0}^{s=\infty} = s \frac{\partial^2}{\partial x^2} \delta(x) \frac{s^{\mu-1} e^{-\frac{s}{2t}}}{t^\mu \Gamma(\mu)} \Bigg|_{s=0}=0, \qquad (\mu >0).
\end{equation*}
where $p(x,0)=\delta(x)$ is the dirac delta function.\\
Thus, 
\begin{align*}
\frac{\partial}{\partial t} q = & \int_0^\infty s \frac{\partial^2}{\partial s^2} \left( p(x,s) \right) Q(s,t)ds + \mu \int_0^\infty \frac{\partial}{\partial s} \left( p(x,s) \right) Q(s,t) ds\\
= & \frac{1}{4} \frac{\partial^4}{\partial x^4} \int_0^\infty s  p(x,s)  Q(s,t)ds + \frac{\mu}{2} \frac{\partial^2}{\partial x^2} \int_0^\infty p(x,s)  Q(s,t) ds.
\end{align*}
In the last calculations we used the fact that the function $p(x,s)$ satisfies the heat equation
\begin{equation*}
\frac{\partial}{\partial s} p(x,s) =\frac{1}{2} \frac{\partial^2}{\partial x^2} p(x,s).
\end{equation*}
Furthermore,
\begin{align*}
\frac{\partial}{\partial t} q =  & \frac{1}{4} \frac{\partial^4}{\partial x^4} \int_0^\infty s  \frac{e^{-\frac{x^2}{2s}}}{\sqrt{2\pi s}}  Q(s,t)ds + \frac{\mu}{2} \frac{\partial^2}{\partial x^2} q(x,t)\\
= & - \frac{1}{4} \frac{\partial^3}{\partial x^3} \int_0^\infty x  \frac{e^{-\frac{x^2}{2s}}}{\sqrt{2\pi s}}  Q(s,t)ds + \frac{\mu}{2} \frac{\partial^2}{\partial x^2} q(x,t)\\
= & -  \frac{1}{4} \frac{\partial^3}{\partial x^3} \left( x q(x,t) \right) + \frac{\mu}{2} \frac{\partial^2}{\partial x^2} q(x,t).
\end{align*}
In order to obtain the formula \eqref{operatorQQ}, further calculations are needed. In particular
\begin{equation*}
\frac{\partial^3}{\partial x^3} \left( x q(x,t) \right) = \left\lbrace 3 \frac{\partial^2}{\partial x^2} + x \frac{\partial^3}{\partial x^3} \right\rbrace q(x,s)
\end{equation*}
and thus
\begin{equation*}
\frac{\partial}{\partial t} q(x,t;\mu) = \left\lbrace \frac{2\mu - 3}{4} \frac{\partial^2}{\partial x^2} - \frac{x}{4} \frac{\partial^3}{\partial x^3} \right\rbrace q(x,t;\mu)
\end{equation*}
This concludes the first proof.
\paragraph{Second proof:}
Let us define the Fourier transform of a function $f$ as follows
\begin{equation}
\mathcal{F}\{ f \} (\beta) = \int_\mathbb{R} e^{-i \beta x} f(x) dx = \left(e^{-i \beta x}, f \right)
\end{equation}
when $\mathcal{F}\{ f \} (\beta)$ exists.\\
Consider now the operator
\begin{equation*}
\mathcal{O}_t(f)= \left\lbrace  (2\mu -3) \frac{\partial^2}{\partial x^2} - x \frac{\partial^3}{\partial x^3} \right\rbrace f,\quad f \in C^\infty_{\downarrow}(\mathbb{R})
\end{equation*}
where $C^\infty_{\downarrow}(\mathbb{R})$ is the space of the rapidly decreasing smooth functions.\\
Evaluate the Fourier transforms
\begin{align*}
\mathcal{F}\left\lbrace  \frac{\partial^2}{\partial x^2} f \right\rbrace (\beta) = \left(e^{-i \beta x}, \frac{\partial^2}{\partial x^2} f \right) = \left( \frac{\partial^2}{\partial x^2} e^{-i \beta x},  f \right) = & - \beta^2 \left(e^{-i \beta x}, f \right) \\
= & -\beta^2 \mathcal{F} \{ f \}(\beta)
\end{align*}
and
\begin{align*}
\mathcal{F}\left\lbrace  x \frac{\partial^3}{\partial x^3} f \right\rbrace (\beta) = & \left(e^{-i \beta x}, x \frac{\partial^3}{\partial x^3} f \right) = \left(x e^{-i \beta x},  \frac{\partial^3}{\partial x^3} f \right) = - \left(\frac{\partial^3}{\partial x^3} x e^{-i \beta x}, f \right) \\
= & 3 \beta^2 \left(e^{-i \beta x}, f \right) + \beta^3 \frac{\partial}{\partial \beta} \left(e^{-i \beta x}, f \right)\\
= & 3 \beta^2 \mathcal{F} \{ f \}(\beta) + \beta^3 \frac{\partial}{\partial \beta} \mathcal{F} \{ f \}(\beta).
\end{align*}
Hence, we have
\begin{equation*}
\mathcal{F}\left\lbrace \mathcal{O}_t(f) \right\rbrace (\beta) = -2\mu \beta^2 \hat{f}(\beta) - \beta^3 \frac{\partial}{\partial \beta} \hat{f}(\beta)
\end{equation*}
where $\hat{f}(\beta)=\mathcal{F}\{ f \}(\beta)$.\\
The Fourier transform of the function \eqref{densityBG1} reads
\begin{equation*}
\hat{q}_t(\beta) = \int_0^\infty e^{-\frac{\beta^2}{2}s} Q(x;t,\mu,1)ds = \left( \frac{1}{1+t\frac{\beta^2}{2}} \right)^\mu = E\left\lbrace e^{- \frac{\beta^2}{2} G_{1,\mu}(t)} \right\rbrace
\end{equation*}
and the time derivative is given by
\begin{equation}
\frac{\partial}{\partial t} \hat{q}_t(\beta) = -\mu \frac{\beta^2}{2} \frac{\hat{q}_t(\beta)}{(1+t\frac{\beta^2}{2})}.
\label{derTimeQQ}
\end{equation}
For the operator \eqref{operatorQQ} we have
\begin{align*}
\mathcal{F}\left\lbrace \frac{1}{4} \mathcal{O}_t(q) \right\rbrace (\beta) = & - \mu \frac{\beta^2}{2} \hat{q}_t(\beta) - \frac{\beta^3}{4} \frac{\partial}{\partial \beta} \hat{q}_t(\beta)\\
= & - \mu \frac{\beta^2}{2} \hat{q}_t(\beta) + \mu \frac{\beta^4}{4} t \frac{1}{(1+t\frac{\beta^2}{2})} \hat{q}_t(\beta)\\
= & - \mu \frac{\beta^2}{2}\hat{q}_t(\beta) \left( 1 -\frac{\frac{\beta^2}{2}t}{(1+t\frac{\beta^2}{2})} \right)\\
= & - \mu \frac{\beta^2}{2} \frac{\hat{q}_t(\beta)}{(1+t\frac{\beta^2}{2})}
\end{align*}
which coincides with \eqref{derTimeQQ}. The proof is completed.

\paragraph{Third proof:}
Evaluate now, the Mellin transform of the function \eqref{densityBG1}
\begin{equation}
\mathcal{M}\left\lbrace q(\cdot \,,t;\mu) \right\rbrace (\eta) = \frac{2^\frac{\eta -1}{2}}{2 \sqrt{\pi}} \Gamma\left( \frac{\eta}{2} \right) \Gamma\left( \frac{\eta -1}{2} +\mu \right) \frac{t^\frac{\eta -1}{2}}{\Gamma(\mu)} = \Psi_t(\eta), \quad \Re\{ \eta \} > 0.
\label{eeT}
\end{equation}
The time derivative of the formula \eqref{eeT} reads
\begin{align*}
\frac{\partial}{\partial t} \Psi_t(\eta) = & \frac{2^\frac{\eta -1}{2}}{2 \sqrt{\pi}} \Gamma\left( \frac{\eta}{2} \right) \Gamma\left( \frac{\eta -1}{2} +\mu \right) \frac{t^\frac{\eta -3}{2}}{\Gamma(\mu)} \left( \frac{\eta -1}{2} \right)\\
= & \frac{2^\frac{\eta -1}{2}}{2 \sqrt{\pi}} \left( \frac{\eta}{2} - 1 \right) \Gamma\left( \frac{\eta -2}{2}\right) \left( \frac{\eta -3}{2} +\mu \right) \Gamma\left( \frac{\eta -3}{2} +\mu \right) \frac{t^\frac{\eta -3}{2}}{\Gamma(\mu)} \left( \frac{\eta -1}{2} \right)\\
= & \frac{1}{4} (\eta -1) (\eta -2) (\eta +2\mu -3) \Psi_t(\eta -2)\\
= & \frac{1}{4} \eta (\eta -1) (\eta -2) \Psi_t(\eta -2) + \frac{1}{4} (\eta -1) (\eta -2) (2\mu -3) \Psi_t(\eta -2)\\
= & -\frac{1}{4} \mathcal{M} \left\lbrace x \frac{\partial^3}{\partial x^3} q(x,t;\mu) \right\rbrace (\eta) + \frac{2\mu -3}{4} \mathcal{M}\left\lbrace \frac{\partial^2}{\partial x^2} q(x,t;\mu) \right\rbrace (\eta).
\end{align*}
The last calculations can be readily verified by exploiting the properties of the Mellin transform \eqref{derM}. This concludes the third proof. $\blacksquare$\\

The moments of the process $B(G_1(t))$, $t>0$ can be evaluated from the Mellin transform \eqref{eeT} as follows
\begin{equation}
E\left\lbrace B(G_1(t)) \right\rbrace^{2k} = 2\Psi_t(2k+1) = \frac{2^k}{\sqrt{\pi}} \Gamma\left( k+\frac{1}{2} \right) \Gamma\left( k + \mu \right) \frac{t^k}{\Gamma(\mu)}, \quad k \in \mathbb{N}.
\end{equation}

The process $B_{\frac{\gamma}{2}}(G_\gamma(t))$, $t>0$ is now examined. In such a process, a fractional Brownian motion $B_{\frac{\gamma}{2}}(t)$, $t>0$ and a generalized Gamma process $G_\gamma(t)$, $t>0$ are involved. The density law of the process $B_{\frac{\gamma}{2}}(G_\gamma(t))$, $t>0$ is given by
\begin{align}
q(x,t)= & \gamma \int_0^\infty \frac{e^{-\frac{x^2}{2s^\gamma}}}{\sqrt{2\pi s^\gamma}} s^{\mu \gamma -1} \frac{e^{-\frac{s^\gamma}{t}}}{t^\mu \Gamma(\mu)}ds \label{ssr}\\
= & 2\frac{|x|^{\mu -\frac{1}{2}} }{ \pi^{\frac{1}{2}} \Gamma(\mu)} \left(\frac{2}{t} \right)^\frac{2\mu +1}{4} K_{\gamma \mu -\frac{\gamma}{2}} \left( |x| \sqrt{\frac{2}{t}} \right), \quad x \in \mathbb{R} \setminus \{0\}, \; t>0, \nonumber
\end{align}
where $\mu >0$ and $\gamma \in (0,2)$.\\ 
The Mellin transform of the function \eqref{ssr} reads
\begin{equation}
\frac{1}{2} E\Big| B_{\frac{\gamma}{2}}(G_\gamma(t)) \Big|^{\eta -1} = \frac{2^\frac{\eta -1}{2}}{2 \sqrt{\pi} \Gamma(\mu)}\Gamma\left( \frac{\eta}{2} \right) \Gamma\left( \frac{\eta -1}{2} + \mu \right) t^\frac{\eta -1}{2}, \quad \Re\{ \eta \} >0 
\end{equation}
We observe that, for $\mu=\frac{1}{2}$, the density law \eqref{ssr} becomes
\begin{equation}
q(x,t)=\frac{1}{\pi} \sqrt{\frac{2}{t}}K_0 \left( |x| \sqrt{\frac{2}{t}} \right), \quad x \in \mathbb{R} \setminus \{ 0 \},\; t>0
\end{equation}
and the parameter $\gamma$ is not relevant. In particular if we choose $\gamma=1$, then we get the result of the Theorem \ref{teQQ} and this result still holds for all $\gamma \in (0,2)$.\\
We also note that 
\begin{equation*}
B_{H}(G_{2H}(t)) \stackrel{i.d.}{=} B(|G_{2H}(t)|^{2H}) \stackrel{i.d.}{=} B(G_1(t)), \quad t>0,\; 0<H<1.
\end{equation*}
Moreover, by \eqref{GdiG}, we can write
\begin{equation}
B(G_1(t)) \stackrel{i.d.}{=} B^1\left( \sqrt{\frac{t}{2}} \right) B^2\left( \sqrt{\frac{t}{2}} \right), \quad t>0
\end{equation}
where $B(t)$, $t>0$, $B^1(t)$, $t>0$ and $B^2(t)$, $t>0$ are independent Brownian motions and $G_1(t)$, $t>0$ is a Gamma process with density law $Q(x;t,\frac{1}{2},1)$.\\

We define and study the $n-$dimensional process
\begin{equation}
I^n_\mu(t) = \left[B^1(G_1(t)), B^2(G_1(t)), \ldots , B^n(G_1(t)) \right]^T , \quad t>0 
\label{ProcmultiBG1}
\end{equation}
where $B^j(t)$, $t>0$, $j=1,2, \ldots , n$ are independent Brownian motions and $G_1(t)$, $t>0$ is a Gamma process. The distribution of the process $I^n_\mu(t)$, $t>0$ is given by
\begin{equation}
q(\textbf{x},t) = \frac{2^{1-\mu}}{\pi^\frac{n}{2}} \frac{\| \textbf{x} \|^{\mu -\frac{n}{2}}}{\Gamma(\mu)}\left( \frac{2}{t}\right)^\frac{2\mu+n}{4} K_{\mu - \frac{n}{2}} \left( \| \textbf{x} \| \sqrt{\frac{2}{t}} \right), \quad \textbf{x} \in \mathbb{R}^n, \; t>0
\label{multiBG1}
\end{equation}
where $\| \textbf{x} \|^2 = \sum_{i=1}^n x_i^2$ is the euclidean norm.\\
In the $2-$dimensional case, we have
\begin{align}
q(x,t;\mu) = & \int_0^\infty \frac{e^{-\frac{x^2}{2s}}}{\sqrt{2\pi s}} \frac{e^{-\frac{y^2}{2s}}}{\sqrt{2\pi s}} \frac{s^{\mu -1} e^{-\frac{s}{t}}}{t^\mu \Gamma(\mu)}ds \label{LLLa2}\\
= & 2^{1-\mu}\frac{(x^2 + y^2)^{\frac{\mu -1}{2}}}{\pi \Gamma(\mu)} \left( \frac{2}{t} \right)^\frac{\mu+1}{2} K_{\mu -1} \left( \sqrt{2 \frac{x^2 + y^2}{t}} \right), \quad (x,y) \in \mathbb{R}^2, \; t>0, \nonumber
\end{align}
where $\mu \neq \frac{1}{2}$. For the case $\mu = \frac{1}{2}$ we refer to the formula \eqref{LLLa}.\\
We provide the following result concerning the process $I^n_\mu(t)$, $t>0$ in the $n-$dimensional case
\begin{te}
The distribution \eqref{multiBG1} of the process $I^n_\mu(t)$, $t>0$ solves the following partial differential equation
\begin{equation}
4 \frac{\partial}{\partial t} q = (2\mu -n)\triangle q - \triangle ( \mathbf{x} \cdot \nabla q ), \quad \mathbf{x} \in \mathbb{R}^n, \; t>0, \; \mu>0.
\label{teMultiDimEq}
\end{equation}
\label{teMultiDim}
\end{te}
\paragraph{Proof:} The distribution \eqref{multiBG1} of the $n-$dimensional process $I_\mu^n(t)$, $t>0$ can be written as follows
\begin{equation*}
q(\textbf{x},t)=\int_0^\infty \prod_{i=1}^n \frac{e^{-\frac{x_i^2}{2s}}}{\sqrt{2\pi s}} \frac{s^{\mu-1} e^{-\frac{s}{t}}}{t^\mu \Gamma(\mu)} ds = \int_0^\infty p(\textbf{x},s) \frac{s^{\mu-1} e^{-\frac{s}{t}}}{t^\mu \Gamma(\mu)} ds , \quad \textbf{x} \in \mathbb{R}^n, \; t>0.
\end{equation*}
At first we point out that
\begin{equation*}
\frac{\partial}{\partial s} p(\textbf{x},s) = \frac{1}{2} \sum_{i=1}^{n} \prod^n_{\begin{subarray}{c} j=1\\ j \neq i \end{subarray}} p(x_j,s) \frac{\partial^2}{\partial x_i^2} p(x_i, s) = \frac{1}{2} \triangle p(\textbf{x}, s)
\end{equation*}
where $\triangle = \sum_{i=1}^n \partial^2_i$ is the Laplacian operator.\\
We are now able to evaluate the time derivative of the function \eqref{multiBG1} which is given by
\begin{align}
\frac{\partial}{\partial t} q(\textbf{x},t) = & [\textrm{by } \eqref{PDEm1}] =  \int_{0}^\infty p(\textbf{x},s) \left\lbrace s \frac{\partial^2}{\partial s^2} Q - (\mu -2) \frac{\partial}{\partial s} Q \right\rbrace ds \label{derivT}\\
= & \int_{0}^\infty p(\textbf{x},s) s \frac{\partial^2}{\partial s^2} Q(s,t) ds - (\mu -2) \int_{0}^\infty p(\textbf{x},t) \frac{\partial}{\partial s} Q(s,t) ds \nonumber \\
= & \mathcal{I}_1 + \mathcal{I}_2 \nonumber
\end{align}
where we used the notation
\begin{equation*}
Q(s,t)=\frac{s^{\mu-1} e^{-\frac{s}{t}}}{t^\mu \Gamma(\mu)}.
\end{equation*}
Let us evaluate the first integral of the time derivative \eqref{derivT}. By performing two integrations by parts we have
\begin{equation*}
\mathcal{I}_1 = p(\textbf{x}, s) s\frac{\partial}{\partial s}Q(s,t) \Bigg|_{s=0}^{s=\infty} - \int_0^\infty  p(\textbf{x},s) \frac{\partial}{\partial s}Q(s,t) ds - \int_0^\infty s \frac{1}{2} \triangle p(\textbf{x},s) \frac{\partial}{\partial s} Q(s,t) ds
\end{equation*}
and
\begin{align*}
\mathcal{I}_1 = & p(\textbf{x},s) s \frac{\partial}{\partial s}Q(s,t) \Bigg|_{s=0}^{s=\infty} - p(\textbf{x},s)Q(s,t) \Bigg|_{s=0}^{s=\infty} \\
 + & \frac{1}{2} \triangle \int_0^\infty p(\textbf{x}, s)Q(s,t)ds - \frac{1}{2} \triangle \int_0^\infty sp(\textbf{x},s) \frac{\partial}{\partial s}Q(s,t) ds
\end{align*}
where we used the fact that $\frac{\partial}{\partial s}p(\textbf{x},s) = \frac{1}{2} \triangle p(\textbf{x},s)$ as we mentioned above.\\
After further integrations by parts, we have
\begin{align*}
\mathcal{I}_1 = & p(\textbf{x},s) s\frac{\partial}{\partial s}Q(s,t) \Bigg|_{s=0}^{s=\infty} - p(\textbf{x},s) Q(s,t) \Bigg|_{s=0}^{s=\infty} + \frac{1}{2} \triangle q(\textbf{x},t)\\
- & \frac{1}{2}\triangle \left\lbrace s p(\textbf{x},s) Q(s,t) \Bigg|_{s=0}^{s=\infty} - \int_0^\infty p(\textbf{x},s)Q(s,t)ds - \int_0^\infty s \frac{\partial}{\partial s}p(\textbf{x},s) Q(s,t)ds \right\rbrace \\
=  & p(\textbf{x},s) s\frac{\partial}{\partial s}Q(s,t) \Bigg|_{s=0}^{s=\infty} - p(\textbf{x},s)Q(s,t) \Bigg|_{s=0}^{s=\infty} - \frac{1}{2} \triangle s p(\textbf{x},s) Q(s,t) \Bigg|_{s=0}^{s=\infty}\\
 + & \triangle q(\textbf{x},t) + \frac{1}{2}\triangle \int_0^\infty s \frac{1}{2}\triangle p(\textbf{x},s) Q(s,t) ds.
\end{align*}
We now observe that
\begin{align*}
p(\textbf{x},s) s\frac{\partial}{\partial s}Q(s,t) \Bigg|_{s=0}^{s=\infty} = & (\mu -1)p(\textbf{x},s) Q(s,t)\Bigg|_{s=0}^{s=\infty} - p(\textbf{x},s) \frac{s}{t}Q(s,t) \Bigg|_{s=0}^{s=\infty} \\
= & (\mu -1)p(\textbf{x},s) Q(s,t) \Bigg|_{s=0}^{s=\infty}
\end{align*}
because
\begin{equation*}
p(\textbf{x},s) \frac{s}{t}Q(s,t) \Bigg|_{s=0}^{s=\infty} = \frac{e^{-\frac{1}{2s} \sum_{i=1}^n x_i^2}}{\sqrt{2\pi s}} \frac{s^\mu e^{-\frac{s}{t}}}{t^{\mu +1} \Gamma(\mu)} \Bigg|_{s=0}^{s=\infty}  = 0, \qquad (\mu >0).
\end{equation*}
Furthermore,
\begin{equation*}
\triangle sp(\textbf{x},s)Q(s,t)\Bigg|_{s=0}^{s=\infty} = - \frac{s^{\mu} e^{-\frac{s}{t}}}{t^\mu \Gamma(\mu)} \triangle \prod_{i=1}^n \delta(x_i) \Bigg|_{s=0} = 0, \qquad (\mu >0)
\end{equation*}
where $p(\textbf{x},0)=\prod_{i=1}^{n}\delta(x_i)$.\\
Hence, the integral $\mathcal{I}_1$ becomes
\begin{equation}
\mathcal{I}_1=(\mu-2)p(\textbf{x},s) Q(s,t) \Bigg|_{s=0}^{s=\infty} + \triangle q(\textbf{x},t) + \frac{1}{4}\triangle \int_0^\infty s \triangle p(\textbf{x},s) Q(s,t) ds.
\label{lastI}
\end{equation}
In order to evaluate the last integral in \eqref{lastI} we point out that 
\begin{align*}
\triangle p(\textbf{x}, s) = & \sum_{i=1}^n \frac{\partial^2}{\partial x_i^2}p(\textbf{x},s) = \sum_{i=1}^n \frac{\partial^2}{\partial x_i^2} \prod_{j=1}^n \frac{e^{-\frac{x_j^2}{2s}}}{\sqrt{2\pi s}} = \sum_{i=1}^n \frac{\partial^2}{\partial x_i^2} \prod_{j=1}^n p(x_j, s)\\
= & \sum_{i=1}^n \prod^n_{\begin{subarray}{c} j=1\\j \neq i \end{subarray}} p(x_j,s) \frac{\partial^2}{\partial x_i^2}p(x_i,s)\\
= & - \sum_{i=1}^n \prod^n_{\begin{subarray}{c} j=1\\j \neq i \end{subarray}} p(x_j,s) \frac{\partial}{\partial x_i} \left( \frac{x_i}{s} p(x_i,s) \right).
\end{align*}
It can be immediately verified that
\begin{equation*}
\frac{\partial^2}{\partial x_i^2} p(x_i, s) = \frac{\partial^2}{\partial x_i^2}\left( \frac{e^{-\frac{x_i^2}{2s}}}{\sqrt{2\pi s}} \right) = -\frac{\partial}{\partial x_i} \left( \frac{x_i}{s} \frac{e^{-\frac{x_i^2}{2s}}}{\sqrt{2\pi s}} \right) = -\frac{\partial}{\partial x_i} \left( \frac{x_i}{s} p(x_i,s) \right).
\end{equation*}
Thus, 
\begin{align*}
\triangle p(\textbf{x}, s) = & - \frac{1}{s} \sum_{i=1}^n \prod^n_{\begin{subarray}{c} j=1\\j \neq i \end{subarray}} p(x_j,s) \left[ p(x_i,s) + x_i \frac{\partial}{\partial x_i} p(x_i,s) \right]\\
= & - \frac{1}{s} \sum_{i=1}^n  \left[ \prod^n_{j=1} p(x_j,s) + x_i \frac{\partial}{\partial x_i} \prod^n_{ j=1} p(x_j,s) \right]\\
= & - \frac{n}{s} p(\textbf{x},s) - \frac{1}{s} \sum_{i=1}^n x_i \frac{\partial}{\partial x_i} p(\textbf{x},s).
\end{align*}
In light of the last result, we can write down
\begin{align*}
\mathcal{I}_1 = & (\mu-2)p(\textbf{x},s) Q(s,t) \Bigg|_{s=0}^{s=\infty} + \triangle q(\textbf{x},t) - \frac{1}{4}\triangle \left( n q(\textbf{x},t) + \sum_{i=1}^n x_i \frac{\partial}{\partial x_i} q(\textbf{x},t) \right) \\
= & (\mu-2)p(\textbf{x},s) Q(s,t) \Bigg|_{s=0}^{s=\infty} + \triangle q(\textbf{x},t) - \frac{n}{4} \triangle q(\textbf{x},t) - \frac{1}{4} \triangle \left( \textbf{x} \cdot \nabla q(\textbf{x},t)  \right).
\end{align*}
The second integral in the formula \eqref{derivT} can be evaluated as follows
\begin{align*}
\mathcal{I}_2 = & - (\mu -2) \int_{0}^\infty p(\textbf{x},t) \frac{\partial}{\partial s} Q(s,t) ds \\
= & - (\mu -2)p(\textbf{x},s)Q(s,t)\Bigg|_{s=0}^{s=\infty} + (\mu -2) \frac{1}{2} \triangle q(\textbf{x},t)
\end{align*}
where we integrated by parts and used the relation $\frac{\partial}{\partial s} p(\textbf{x},s) = \frac{1}{2} \triangle p(\textbf{x},s)$.\\
Finally, by combining all the previous results, we have
\begin{equation*}
\mathcal{I}_1 + \mathcal{I}_2 = \frac{2\mu -n}{4} \triangle q(\textbf{x},t) - \frac{1}{4} \triangle \left( \textbf{x} \cdot \nabla q(\textbf{x},t)  \right)
\end{equation*}
and this concludes the proof. $\blacksquare$\\

\begin{os}
\normalfont
We focus now on the Theorem \ref{teMultiDim}, in particular on the case where the shape parameter $\mu=\frac{n}{2}$. From the formulae \eqref{propertiesG} and \eqref{propSomma} we know that
\begin{equation}
\sum_{i=1}^n |B^i(t)|^2 \stackrel{i.d.}{=} \sum_{i=1}^n G^i_{1,\frac{1}{2}}(2t) \stackrel{i.d.}{=} G_{1,\frac{n}{2}}(2t), \quad t>0.
\end{equation}
Hence, we can write down
\begin{equation}
B(G_{1,\frac{n}{2}}(2t)) \stackrel{i.d.}{=} B(BSQ^0_n(t)), \quad t>0
\end{equation}
where $BSQ^0_n(t)$, $t>0$ is the $n-$dimensional squared Bessel process starting from zero.\\
Squared Bessel processes are Markov processes, and their transition densities are known explicitly. For the $\delta-$dimensional $BSQ^x_\delta(t)$, $t>0$ starting from $x \geq 0$ the transition density is given by
\begin{equation}
\label{funBSQ}
q(y,t;x) = \frac{1}{2t} \left( \frac{y}{x} \right)^\frac{\delta -2}{4} e^{-\frac{x+y}{2t}} I_{\frac{\delta -2}{2}} \left( \frac{\sqrt{xy}}{t} \right), \quad y>0, \; t>0, \; \delta >0.
\end{equation}
$I_\nu(z)$ is the Bessel modified function of the first kind (see formula \eqref{I0Bessel1}).\\
The function \eqref{funBSQ} reduces to the Gamma distribution when the starting point is $x=0$, and we obtain
\begin{equation}
q(y,t;0)=\frac{y^{\frac{\delta}{2} -1} e^{-\frac{y}{2t}}}{(2t)^\frac{\delta}{2} \Gamma\left( \frac{\delta}{2} \right)}, \quad y>0, \; t>0, \; \delta>0.
\end{equation}
Now, for $\mu=\frac{n}{2}$ in the p.d.e. \eqref{teMultiDimEq}, we have a new p.d.e. which writes
\begin{equation}
4 \frac{\partial}{\partial t} q = - \triangle (\textbf{x} \cdot \nabla q), \quad \textbf{x} \in \mathbb{R}^n, \; t>0
\label{newPDE}
\end{equation}
and, for $q=q(x,2t)$, we have
\begin{equation}
\frac{\partial}{\partial t} q = - \frac{1}{2} \triangle (\textbf{x} \cdot \nabla q), \quad \textbf{x} \in \mathbb{R}^n, \; t>0
\label{newPDE2}
\end{equation}
The distribution
\begin{equation}
q(\textbf{x}, t) = \frac{2}{(\pi t)^\frac{n}{2} \Gamma\left( \frac{n}{2} \right)} K_0 \left( \frac{\| \textbf{x} \|}{ \sqrt{t}} \right), \quad \textbf{x} \in \mathbb{R}^{n}, \; t>0
\label{sDi}
\end{equation}
which comes from the formula \eqref{multiBG1}, satisfies the p.d.e \eqref{newPDE2}.\\
Finally the p.d.e. \eqref{newPDE2} is the governing equation of the $n-$dimensional process
\begin{equation}
I_\frac{n}{2}^n(2t) = \left( B^1(BSQ^0_n(t), B^2(BSQ^0_n(t), \ldots , B^n(BSQ^0_n(t)) \right)^T, \quad t>0
\end{equation}
which possesses the distribution \eqref{sDi}.\\
For $n=1$ we have
\begin{equation}
I^1_\frac{1}{2}(2t) = B(G_{1,\frac{1}{2}}(2t)) \stackrel{i.d.}{=} B^1(| B^2(t)|^2), \quad t>0
\end{equation}
and the governing equation becomes 
\begin{equation}
\frac{\partial}{\partial t} q = -\frac{1}{2} \frac{\partial^2}{\partial x^2} \left( x \frac{\partial}{\partial x} q \right), \quad x \in \mathbb{R}, \, t>0.
\label{Wpde}
\end{equation}
\end{os}

\begin{os}
\normalfont
Consider the process
\begin{equation}
\| I^n_\mu (t)\| = \sqrt{\sum_{j=1}^{n} \bigg| B^j(G_{1, \mu}(t)) \bigg|^2 }, \quad t>0
\label{procZZ}
\end{equation}
where $B^j(t)$, $t>0$ , $j=1,2,\ldots , n$ are independent Brownian motions starting from zero and $G_{1,\mu}(t)$, $t>0$ is a Gamma process. The process \eqref{procZZ} is the euclidean norm of the process in \eqref{ProcmultiBG1}, $I^n_\mu(t)$, $t>0$.\\
We can find the distribution of the process $\| I^n_\mu (t)\|$, $t>0$ by observing that
\begin{align*}
\| I^n_\mu (t)\| \stackrel{i.d.}{=} & [ \textrm{by } \eqref{propertiesG} \textrm{ and } \eqref{ppppA} ] \stackrel{i.d.}{=} \sqrt{\sum_{j=1}^{n} \, 2 \,G^j_{1, \frac{1}{2}}(G_{1, \mu}(t)) }, \quad t>0\\
\stackrel{i.d.}{=} & [ \textrm{by } \eqref{prodGn} ] \stackrel{i.d.}{=} \sqrt{2\, \sum_{j=1}^{n} G^j_{1, \frac{1}{2}}(t^\frac{1}{2})G_{1, \mu}(t^\frac{1}{2}) }, \quad t>0 \\
\stackrel{i.d.}{=} & \sqrt{ 2\, G_{1, \mu}(t^\frac{1}{2}) \sum_{j=1}^{n} G^j_{1, \frac{1}{2}}(t^\frac{1}{2}) }, \quad t>0\\
\stackrel{i.d.}{=} & [ \textrm{by } \eqref{propSomma} ] \stackrel{i.d.}{=} \sqrt{2\, G_{1, \mu}(t^\frac{1}{2}) G_{1, \frac{n}{2}}(t^\frac{1}{2}) }, \quad t>0\\
\stackrel{i.d.}{=}& [ \textrm{by } \eqref{prodGn} \textrm{ and } \eqref{ppppA} ] \stackrel{i.d.}{=} \sqrt{ G_{1, \frac{n}{2}}(2\,G_{1, \mu}(t) ) }\\
\stackrel{i.d.}{=} & [ \textrm{by } \eqref{selfGcomp} \textrm{ and } \eqref{propertiesG} ] \stackrel{i.d.}{=} G_{2, \frac{n}{2}}(2\,G_{1, \mu}(t) ) , \quad t>0.
\end{align*}
Therefore, the distribution of the process $\| I^n_\mu (t)\|$, $t>0$ follows from the formula \eqref{specialGG}. The process $\| I^n_\mu (t)\| \stackrel{i.d.}{=}\sqrt{2\, G^1_{1, \mu}(t^\frac{1}{2}) G^2_{1, \frac{n}{2}}(t^\frac{1}{2}) }$, $t>0$ can be also written as 
\begin{align*}
\| I^n_\mu (t)\| \stackrel{i.d.}{=} & \sqrt{2\, G^1_{1, \mu}(G^2_{1, \frac{n}{2}}(t)) }, \quad t>0\\
\stackrel{i.d.}{=} & [ \textrm{by } \eqref{selfGcomp} \textrm{ and } \eqref{propertiesG} ] \stackrel{i.d.}{=} \sqrt{2} \,G^1_{2, \mu}(G^2_{1, \frac{n}{2}}(t)) , \quad t>0\\
\stackrel{i.d.}{=} & [ \textrm{by } \eqref{ppppA} ] \stackrel{i.d.}{=} G^1_{2, \mu}(2G^2_{1, \frac{n}{2}}(t)) , \quad t>0\\
\stackrel{i.d.}{=} & G^1_{2, \mu}(2BSQ^0_{n}(t/2)) , \quad t>0.
\end{align*}
Hence, it follows that 
\begin{equation}
\|I^n_\frac{1}{2}(2t) \| \stackrel{i.d.}{=} \big| B(BSQ^0_{n}(t)) \big|, \quad t>0.
\label{JJla}
\end{equation}
The process \eqref{JJla} can be written as follows
\begin{align*}
\|I^n_\frac{1}{2}(2t) \| = & \| \left( B^1(|B(t)|^2), B^2(|B(t)|^2), \ldots , B^n(|B(t)|^2 ) \right)^T \|, \quad t>0\\
=& \| \mathbf{B}( | B(t) |^2) \|, \quad t>0
\end{align*}
and thus,
\begin{equation}
\| \mathbf{B}( | B(t) |^2) \| \stackrel{i.d.}{=} \big| B( \| \mathbf{B}(t) \|^2) \big|, \quad t>0
\end{equation}
where $\mathbf{B}(t)=\left( B^1(t), B^2(t), \ldots , B^n(t) \right)^T$, $t>0$. 
\end{os}

\begin{os}
\normalfont
We consider the composition $\tilde{G}^1_{\gamma, \frac{1}{2}}(\tilde{G}^2_{-\gamma, \frac{\nu}{2}}(t))$, $t>0$ with distribution \eqref{tdist1} examined in the previous section,  where two independent generalized Gamma processes are involved. The process $T(t)=\tilde{G}^1_{2, \frac{1}{2}}(\tilde{G}^2_{-2, \frac{\nu}{2}}(t^\frac{1}{2}))$, $t>0$, for $\gamma=2$, has a distribution given by
\begin{equation}
q(x,t)=2 \frac{t^\frac{\nu}{2}}{(x^2 + t)^\frac{\nu +1}{2}} \frac{\Gamma\left( \frac{\nu+1}{2} \right)}{ \sqrt{\pi} \Gamma\left( \frac{\nu}{2} \right)}, \quad x \geq 0, \; t>0, \; \nu >0.
\end{equation}
Furthermore, by keeping in mind the property $\tilde{G}_{\gamma, \mu}(t^\frac{1}{\gamma}) \stackrel{i.d.}{=} G_{\gamma, \mu}(t)$, $t>0$, $\gamma \neq 0$ we have
\begin{align}
T(t) \stackrel{i.d.}{=} &  \tilde{G}^1_{2, \frac{1}{2}}( G^2_{-2, \frac{\nu}{2}}(t^{-1})), \quad t>0\\
\stackrel{i.d.}{=} & G^1_{2, \frac{1}{2}}(| G^2_{-2, \frac{\nu}{2}}(t^{-1}) |^2), \quad t>0\nonumber \\
\stackrel{i.d.}{=} & [ \textrm{by } \eqref{propertiesG} ] \stackrel{i.d.}{=} G^1_{2, \frac{1}{2}}( G^2_{-1, \frac{\nu}{2}}(t^{-1}) ), \quad t>0 \nonumber\\
\stackrel{i.d.}{=} & [ \textrm{by } \eqref{ppppA} ] \stackrel{i.d.}{=} G^1_{2, \frac{1}{2}}(2 G^2_{-1, \frac{\nu}{2}}(2t^{-1}) ), \quad t>0 \nonumber\\
\stackrel{i.d.}{=} & \Bigg| B\left( \frac{1}{G_{1, \frac{\nu}{2}}\left( \frac{2}{t}\right)} \right) \Bigg|, \quad t>0, \nonumber\\
\stackrel{i.d.}{=} & \Bigg| B\left( \frac{1}{BSQ^0_{\nu}\left( \frac{1}{t}\right)} \right) \Bigg|, \quad t>0.
\end{align}
where the process $BSQ^0_{\nu}(t)$, $t>0$ is the $\nu-$dimensional squared Bessel process starting from the origin.\\
Furthermore, from the property 
\begin{equation}
\tilde{G}^1_{\gamma, \mu_1}(\tilde{G}^2_{-\gamma, \mu_2}(t)) \stackrel{i.d.}{=} \tilde{G}^1_{\gamma, \mu_1}(t^\frac{1}{2}) \, \tilde{G}^2_{-\gamma, \mu_2}(t^\frac{1}{2}), \quad t>0
\end{equation}
we are able to write down
\begin{align}
T(t^2) \stackrel{i.d.}{=} & \tilde{G}^1_{2, \frac{1}{2}}(t^\frac{1}{2}) \, \tilde{G}^2_{-2, \frac{\nu}{2}}(t^\frac{1}{2}), \quad t>0\\
\stackrel{i.d.}{=} & G^1_{2, \frac{1}{2}}(t) \, G^2_{-2, \frac{\nu}{2}}(t^{-1}), \quad t>0\nonumber \\
\stackrel{i.d.}{=} & [ \textrm{by } \eqref{invaaaa} ] \stackrel{i.d.}{=} \frac{G^1_{2, \frac{1}{2}}(t)}{G^2_{2, \frac{\nu}{2}}\left( \frac{1}{t} \right)}\stackrel{i.d.}{=} [ \textrm{by } \eqref{ppppA} ] \stackrel{i.d.}{=} \frac{G^1_{2, \frac{1}{2}}(2t)}{G^2_{2, \frac{\nu}{2}}\left(\frac{2}{t}\right)}, \quad t>0\nonumber \\
\stackrel{i.d.}{=} & \Bigg| \frac{B(t)}{BS^0_\nu \left(\frac{1}{t} \right)} \Bigg|, \quad t>0 \nonumber
\end{align}
where $BS^0_\nu (t)$, $t>0$ is the $\nu-$dimensional Bessel process starting from zero.\\
We observe that (by the symmetry of the distributions)
\begin{equation}
B\left( \frac{1}{BSQ^0_{\nu}\left( \frac{1}{t^2}\right)} \right) = B\left( \Bigg| \frac{1}{BS^0_{\nu}\left( \frac{1}{t^2}\right)} \Bigg|^2 \right) \stackrel{i.d.}{=} \frac{B(t)}{BS^0_\nu \left(\frac{1}{t} \right)}, \quad t>0.
\end{equation}
We point out that, for each $\nu>0$, the random variable $T(\nu)$ possesses the density law $q(x,\nu)$ which is the distribution of a Student's random variable with $\nu>0$ degrees of freedom.
\end{os}

We have already studied the process $B(G_1(t))$, $t>0$ where $B(t)$, $t>0$ is a standard Brownian motion and $G_1(t)$, $t>0$ is a Gamma process with distribution depending on the parameter $\mu>0$. A particular case where $\mu=\frac{1}{2}$ is now examined.\\ 
Consider the process $B^1_{H_1}(|B^2_{H_2}(t)|^\frac{1}{H_1})$, $t>0$  where $B^j_{H_j}(t)$, $t>0$, $j=1,2$ are independent fractional Brownian motions with covariance function given by
\begin{equation}
E\left\lbrace B^j_{H_j}(t) B^j_{H_j}(s)  \right\rbrace = \frac{1}{2} \left( |t|^{2H_j} + |s|^{2H_j} - |t-s|^{2H_j} \right), \quad j=1,2
\end{equation}
for $0< H_1$, $H_2<1$.\\
The distribution of the process $B^1_{H_1}(|B^2_{H_2}(t)|^\frac{1}{H_1})$, $t>0$ reads
\begin{equation}
q(x,t)=2 \int_0^\infty \frac{e^{-\frac{x^2}{2s^2}}}{\sqrt{2\pi s^2}} \frac{e^{-\frac{s^2}{2t^{2H}}}}{\sqrt{2\pi t^{2H}}} ds = \frac{1}{\pi t^H} K_0\left( \frac{|x|}{t^H} \right), \quad x \in \mathbb{R} \setminus \{ 0 \}, \; t>0
\label{densityK0}
\end{equation}
where, for the sake of simplicity, we set $H_2=H$.\\ 
The distribution \eqref{densityK0} comes from  the formula \eqref{densityUNOK} with $\gamma=2$ and $\mu=\frac{1}{2}$ or the formula \eqref{densityBG1} with $\mu=\frac{1}{2}$ by considering the symmetry of the distribution.\\

We give the following result (see, for more details \cite{DO})
\begin{te}
The process $B^1_{H_1}(|B^2_{H_2}(t)|^\frac{1}{H_1})$, $t>0$ has a density law  \eqref{densityK0} which solves the partial differential equation
\begin{equation}
\frac{\partial}{\partial t} q=-H_2t^{2H_2-1} \left( 2\frac{\partial^2}{\partial x^2} + x \frac{\partial^3}{\partial x^3} \right) q
\label{W2pde}
\end{equation}
for all $x\neq 0$, $t>0$ and $H_2 \in (0,1)$.
\end{te}
\paragraph{Proof:} 
From the Theorem \ref{teQQ} we know that the distribution of the process $B(G_1(t))$, $t>0$ is solution to the p.d.e.
\begin{equation*}
4\frac{\partial}{\partial t}Q = (2\mu -3) \frac{\partial^2}{\partial x^2}Q - x \frac{\partial^3}{\partial x^3} Q, \quad x \in \mathbb{R},\; t>0, \; \mu>0.
\end{equation*}
In the special case where $\mu=\frac{1}{2}$ we have
\begin{equation*}
4\frac{\partial}{\partial t}Q = -2 \frac{\partial^2}{\partial x^2}Q - x \frac{\partial^3}{\partial x^3} Q, \quad x \in \mathbb{R},\; t>0.
\end{equation*}
We also know that $B^1_{H_1}(|B^2_{H_2}(t)|^\frac{1}{H_1}) \stackrel{i.d.}{=} B(G_1(2t^{2H}))$ where $G_1$ has a distribution $Q(x;2t^{2H}, \frac{1}{2},1)$ and thus
\begin{equation*}
\frac{\partial}{\partial t}q =\frac{d}{d t}\left( 2t^{2H} \right) \frac{\partial}{\partial z} Q \bigg|_{z=2t^{2H}} = -Ht^{2H-1} \left( 2\frac{\partial^2}{\partial x^2} + x \frac{\partial^3}{\partial x^3} \right) q, \quad x \in \mathbb{R}, \; t>0.
\end{equation*}
Moreover, from \eqref{K0in0}, we have to observe that
\begin{equation*}
\lim_{x \to 0^+} q(x,t) = \infty.
\end{equation*} 
The proof is completed. $\blacksquare$\\

For the process $B^1_{H}(|B^2_{H}(t)|^\frac{1}{H}) \stackrel{i.d.}{=} B^1(|B^2_{H}(t)|^2)$, $t>0$ we prove that
\begin{equation}
E\left\lbrace B^1_{H}(|B^2_{H}(t)|^\frac{1}{H}) \right\rbrace^{\eta -1} = E\left\lbrace B^1_{\frac{H}{2}}(t^\frac{1}{2}) \right\rbrace^{\eta -1} E\left\lbrace B^2_{\frac{H}{2}}(t^\frac{1}{2}) \right\rbrace^{\eta -1}
\end{equation}
and thus we can write down
\begin{equation}
B^1_{H}(|B^2_{H}(t)|^\frac{1}{H}) \stackrel{i.d.}{=} B^1_{\frac{H}{2}}(t^\frac{1}{2}) B^2_{\frac{H}{2}}(t^\frac{1}{2}), \quad t>0.
\label{LLkk}
\end{equation}
We observe that
\begin{align*}
\left\lbrace B^1(|B^2_{H}(t)|^2)\, B^1(|B^2_{H}(s)|^2) \right\rbrace \stackrel{i.d.}{=} & \left\lbrace B^1_\frac{H}{2}(t)\,B^2_\frac{H}{2}(t)\,B^1_\frac{H}{2}(s)\,B^2_\frac{H}{2}(s) \right\rbrace 
\end{align*}
(by \eqref{LLkk}) and 
\begin{align}
E\left\lbrace B^1_\frac{H}{2}(t)\,B^2_\frac{H}{2}(t)\,B^1_\frac{H}{2}(s)\,B^2_\frac{H}{2}(s) \right\rbrace = &  E\left\lbrace B^1_\frac{H}{2}(t)\, B^1_\frac{H}{2}(s) \right\rbrace \, E\left\lbrace B^2_\frac{H}{2}(t)\,B^2_\frac{H}{2}(s) \right\rbrace \nonumber \\
= & \frac{1}{4} \left( |t|^{H} + |s|^{H} - |t-s|^{H} \right)^2.
\label{CovBBH}
\end{align}
by the independence of the processes $B^1_{H}(t)$, $t>0$ and $B^2_{H}(t)$, $t>0$.\\
The covariance function \eqref{CovBBH} can be only positive.\\ 

We have already seen in the formula \eqref{Wpde} a special case of the p.d.e. \eqref{W2pde} where the subordination is made with respect to the fractional Brownian motion insteed of the Brownian motion.\\

In a general setting we can state the following result (see, for more details \cite{DO}).\\
For the process 
\begin{equation}
I^{n-1}_F (t)= B_H^{1} (| B_H^{2}( \ldots | B^{n}_{H}(t) |^{\frac{1}{H}} \ldots ) |^\frac{1}{H}), \quad  t>0, \; H \in (0,1)  
\end{equation}
with $B^j_H$, $j=1,2,\ldots , n$ independent fractional Brownian motions, the following equalities in distribution hold
\begin{equation}
I^{n-1}_F (t) \stackrel{i.d.}{=} B^{1} (| B^{2}( \ldots | B^{n}_{H}(t) |^2 \ldots ) |^2) \stackrel{i.d.}{=} \prod_{i=1}^{n} B^i_{\frac{H}{n}} (t), \quad t>0, \; H \in (0,1).
\label{GhGhG}
\end{equation}
\label{teEXTEND}
The proof of the first equality follows from the property $B_H(t) \stackrel{i.d.}{=} B(t^{2H})$ for each $t$. The proof of the second equality can be given by evaluating the Mellin transform of the density of $I^{n-1}_F(t)$, $t>0$ which reads
\begin{align*}
E \left\lbrace I^{n-1}_F (t) \right\rbrace^{\eta -1} = & \int_{0}^{\infty} x^{\eta -1} 2^{n-1} \int_{\mathbb{R}^{n-1}} \frac{e^{-\frac{x^{2}}{2s_1^2}}}{\sqrt{2 \pi s_1^2}} \frac{e^{-\frac{s^{2}_1}{2s_2^2}}}{\sqrt{2 \pi s_2^2}} \ldots \frac{e^{-\frac{s^{2}_{n-1}}{2t^{2H}}}}{\sqrt{2 \pi t^{2H}}}ds_1 \ldots ds_{n-1}\\
= & \left[ \frac{2^{\frac{\eta}{2}} \Gamma\left( \frac{\eta}{2} \right)}{\sqrt{2\pi}} \right]^n t^{H(\eta -1)}.\\
= & \prod_{i=1}^n E \left\lbrace B^i_{\frac{H}{n}} (t)  \right\rbrace^{\eta -1}.
\end{align*}
Such a transform is in line with the Mellin convolution machinery.\\
Also, we remind that
\begin{equation*}
|B_H(t)|^2 \stackrel{i.d.}{=} |G_{2, \frac{1}{2}}(g(t))|^2 \stackrel{i.d.}{=} G_{1,\frac{1}{2}}(g(t)), \quad t>0
\end{equation*}
where $g(t)=2t^{2H}$.\\
Hence, we have
\begin{align*}
I^{n-1}_F(t) \stackrel{i.d.}{=} & B\left( G^1_{1,\frac{1}{2}}(G^2_{1,\frac{1}{2}}( \ldots G^{n-1}_{1,\frac{1}{2}}(g(t)) \ldots) \right), \quad t>0\\
\stackrel{i.d.}{=} & G^n_{2,\frac{1}{2}}\left( G^1_{1,\frac{1}{2}}(G^2_{1,\frac{1}{2}}( \ldots G^{n-1}_{1,\frac{1}{2}}(g(t)) \ldots) \right), \quad t>0\\
\stackrel{i.d.}{=} & [ \textrm{by the Theorem } \ref{teP1} ]\\
\stackrel{i.d.}{=} & \prod_{j=1}^n G^j_{2,\frac{1}{2}} \left(  |g(t)|^\frac{1}{n} \right), \quad t>0\\
\stackrel{i.d.}{=} & \prod_{j=1}^n B^j_{H} \left(  t^\frac{1}{n} \right), \quad t>0\\ 
\stackrel{i.d.}{=} & \prod_{j=1}^n B^j_\frac{H}{n} \left( t \right), \quad t>0.
\end{align*}
The last two steps follow from the property $B_H(t) \stackrel{i.d.}{=} B(t^{2H})$, for each $t>0$.\\ 

From the covariance function \eqref{CovBBH}, by considering the Theorem \ref{teEXTEND} we obtain
\begin{align*}
\left\lbrace I_F^{n-1}(t) I_F^{n-1}(s) \right\rbrace \stackrel{i.d.}{=} & [\textrm{by } \eqref{GhGhG}] = E\left\lbrace \prod_{i=1}^n B^i_\frac{H}{n}(t) \prod_{i=1}^n B^i_\frac{H}{n}(s) \right\rbrace
\end{align*}
and
\begin{align}
E\left\lbrace \prod_{i=1}^n B^i_\frac{H}{n}(t) \prod_{i=1}^n B^i_\frac{H}{n}(s) \right\rbrace = & [\textrm{by } \eqref{GhGhG}] = E\left\lbrace \prod_{i=1}^n B^i_\frac{H}{n}(t) \prod_{i=1}^n B^i_\frac{H}{n}(s) \right\rbrace \nonumber \\
= & \prod_{i=1}^n E\left\lbrace B^i_\frac{H}{n}(t) \, B^i_\frac{H}{n}(s) \right\rbrace \nonumber \\
= & \frac{1}{2^n} \left( |t|^\frac{2H}{n} + |s|^\frac{2H}{n} - |t-s|^\frac{2H}{n} \right)^n. \nonumber
\end{align}

We observe that the $n-$dimensional process 
\begin{equation}
\left( B^1(|B(t)|^2), B^2(|B(t)|^2), \ldots , B^n(|B(t)|^2) \right)^T, \quad t>0
\end{equation} 
is a special case of the process \eqref{ProcmultiBG1}, in particular $I^2_\frac{1}{2}(2t)$, $t>0$. Such a process possesses a distribution given by the formula \eqref{multiBG1} with $\mu=\frac{1}{2}$ which satisfies the partial differential equation
\begin{equation}
2 \frac{\partial}{\partial t} q = (1 -n)\triangle q - \triangle ( \textbf{x} \cdot \nabla q ), \quad \textbf{x} \in \mathbb{R}^n, \; t>0.
\end{equation}

In the $2-$dimensional case we have the distribution \eqref{LLLa2} which becomes 
\begin{align}
q(\textbf{x},t) = & \frac{2^\frac{1}{4}}{\pi \sqrt{\pi}} \frac{\| \textbf{x} \|^{-\frac{1}{2}}}{t^\frac{3}{4}} K_{-\frac{1}{2}} \left( \|\textbf{x} \| \sqrt{\frac{2}{t}}\right) \label{LLLa}\\
= & \frac{2^\frac{1}{4}}{\pi \sqrt{\pi}} \frac{\| \textbf{x} \|^{-\frac{1}{2}}}{t^\frac{3}{4}} \sqrt{\frac{\pi}{2}} \left( \frac{\sqrt{t}}{\| \textbf{x} \| \sqrt{2}} \right)^\frac{1}{2} e^{-\|\textbf{x} \| \sqrt{\frac{2}{t}}} \nonumber \\
= & \frac{1}{\pi \sqrt{t}} \frac{1}{\| \textbf{x} \|} e^{-\|\textbf{x} \| \sqrt{\frac{2}{t}}} \nonumber
\end{align}
where we used the property $K_{-\frac{1}{2}}(z)=K_{\frac{1}{2}}(z)=\sqrt{\frac{\pi}{2z}}e^{-z}$ (see p. 925 Gradshteyn and Ryzhik \cite{GR}).\\
Therefore, it follows that
\begin{align*}
p(\textbf{x},t) = & q(\textbf{x},2t)\\
= & \frac{1}{2\pi \sqrt{t}} \frac{e^{-\sqrt{\frac{(x_1^2 +x_2^2 )}{t}}}}{\sqrt{(x_1^2 + x_2^2) }}, \quad (x_1, x_2) \in \mathbb{R}^2 \setminus \{ 0, 0 \}, \; t>0.
\end{align*}
is the distribution of the process $\|I^2_\frac{1}{2}(2t)\|= \sqrt{|B^1(|B(t)|^2)|^2 + |B^2(|B(t)|^2)|^2}$, $t>0$.\\

We give now some connections between the iterated fractional Brownian motion and the composition of the generalized Gamma processes.\\ 
Consider the process 
\begin{equation}
P_n(t)=B\left( \prod_{j=1}^{n} G^j_{\gamma_j}(2t) \right) , \quad t>0
\end{equation}
where $B(t)$, $t>0$ is a standard Brownian motion and $G^j_{\gamma_j}(t)$, $t>0$, $j=1,2,\ldots ,n$ are independent generalized Gamma processes. The process $P_n(t)$, $t>0$ possesses a density law given  by
\begin{equation}
q(x,t)=\int_{0}^\infty \frac{e^{-\frac{x^2}{2s}}}{\sqrt{2\pi s}} Q^n(s;2t,\mu,\gamma_1, \gamma_2, \ldots , \gamma_n) ds, \quad x \in \mathbb{R}, \; t>0
\label{dens1}
\end{equation}
where $Q^n(s;2t,\mu,\gamma_1, \gamma_2, \ldots , \gamma_n)$ is the distribution of the product $\prod_{j=1}^{n} G^j_{\gamma_j}(2t)$, $t>0$.\\
Consider now $\gamma_j=\frac{1}{H_j}$ for $j=1,2,\ldots , n$. The Mellin transform of the formula \eqref{dens1} is given by
\begin{equation}
\mathcal{M}\left\lbrace q(x,t) \right\rbrace (\eta) = \frac{2^\frac{\eta -1}{2}}{\sqrt{\pi}} \Gamma\left( \frac{\eta}{2} \right) \int_0^\infty s^{\eta -1} Q^n(s; 2t, \mu, \gamma_1, \gamma_2, \ldots, \gamma_n) ds.
\end{equation}
The last integral is the Mellin transform of the process  \eqref{IGn} and is given by the formula \eqref{mellinGn} which coincides, for $\mu=\frac{1}{2}$, with the Mellin transform of the distribution of the iterated fractional Brownian motion. Thus, we can write down the following equality in distribution
\begin{equation}
B\left( \prod_{j=1}^{n} G^j_\frac{1}{H_j}(2t) \right) \stackrel{i.d.}{=} B^1_{H_1}(|B^2_{H_2}(\ldots |B^{n+1}_{H_{n+1}}(t) | \ldots)|), \quad t>0.
\label{H1}
\end{equation}

Furthermore, for $\mu=1/2$ and $\gamma_j=1$ for $j=1,2,\ldots , n$ (or equivalently $H_j=1$ for $j=1,2,\ldots , n$ in the process \eqref{H1}) we have
\begin{align*}
B\left( \prod_{j=1}^{n} G^j_{1}(2t) \right) \stackrel{i.d.}{=} & B^1_{H_1}(|B^2_{H_2}(\ldots |B_{H_{n+1}}^{n+1}(t)|^{1/H_n} \ldots)|^{1/H_1}), \quad t>0.\\
\stackrel{i.d.}{=} &  B^1(|B^2(|\ldots |B^{n+1}_H(t)|^2 \ldots |)|^2), \quad t>0
\end{align*}
where $H_{n+1}=H \in (0,1)$.

\end{document}